\documentclass[10pt]{article}

\usepackage{a4wide}
\usepackage{amssymb}
\usepackage{amsfonts}
\usepackage{amsmath}
\input xy
\xyoption{arrow} \xyoption{matrix}

\date{}

\newtheorem{proposition}{Proposition}[section]
\newtheorem{theorem}[proposition]{Theorem}
\newtheorem{lemma}[proposition]{Lemma}

\newtheorem{corollary}[proposition]{Corollary}

\def\GK{{\rm  GK}\,}

\def\Hom{{\rm Hom}}
\def\der{\partial }

\def\nFM0{{\nu }_{F,M_0}}
\def\nFN0{{\nu }_{F,N_0}}
\def\nGN0{{\nu }_{G,N_0}}

\def\N0{ {\bf N}_0 }

\def\t{\otimes}
\def\g{\gamma}
\def\v{\varphi}
\def\ra{\rightarrow}

\def\lra{\leftrightarrow}
\def\Xpm{X^{\pm }}

\def\s{\sigma}
\def\Z{\mathbb{Z}}

\def\l1{{\lambda}_1}

\def\a{\alpha}
\def\a0{ {\alpha }_0}
\def\a1{ {\alpha }_1}

\def\l{\lambda}
\def\o{\omega}

\def\nFGM0{{\nu }_{F,G,M_0}}

\def\nFN0{{\nu}_{F,N_0}}

\def\sm{{\sigma}^m}

\def\sm1{{\sigma}^{-1}}

\def\smtp1{{\sigma}^{-t+1}}

\def\o{\omega }
\def\S1{S^{-1}}

\def\Xpm1{X^{\pm 1}_1}

\def\sPM1{{\sigma }^{\pm 1}}
\def\sMP1{{\sigma }^{\mp 1 }}

\def\d{\delta}

\def\di{{\rm d.ind}}

\def\L{\Lambda}

\def\G{\Gamma}

\def\Ytm1{Y^{t-1}}
\def\Yim1{Y^{i-1}}

\def\CK{{\cal K}}

\def\CN{{\cal N}}

\def\CH{{\cal H}}

\def\CZ{{\cal Z}}

\def\Aut{{\rm Aut}}

\def\dim{{\rm dim }}

\def\ker{ {\rm ker } }

\def\CJ{ {\cal J}}

\def\D{ \Delta }

\def\SL2Z{ {\rm SL}_2({\bf Z}) }

\def\CZ{ {\cal Z}}

\def\th{ \theta }

\def\Gp1{ G^{1 , 1 } }
\def\P11{ P^{-1 , 1 } }
\def\Pp1{ P^{1 , 1 } }

\def\th{\theta}

\def\nCLsr{{}^\nu\kern-2pt {\cal L}^{\sigma , \rho  }}
\def\nP{{}^\nu \kern-2pt P}
\def\nL{{}^\nu\kern-2pt L}
\def\nLL{{}^\nu\kern-2pt \Lambda}
\def\nPsr{{}^\nu\kern-2pt P^{\sigma , \rho  }}
\def\nLsr{{}^\nu\kern-2pt L^{\sigma , \rho  }}
\def\nuCL{{}^\nu\kern-2pt  {\cal L}}
\def\nCLsr{{}^\nu\kern-2pt {\cal L}^{\sigma , \rho  }}
\def\nCL1m{{}^\nu\kern-2pt {\cal L}^{-1 , 1  }}
\def\x1nu{x^\frac{1}{\nu}}
\def\xm1nu{x^{-\frac{1}{\nu}}}

\def\CN{{\cal N}}
\def\ra{\rightarrow }

\def\CI{{\cal I}}

\def\coker{{\rm coker}}

\def\CC{ {\cal C}}

\def\CH{ {\cal H}}
\def\CP{ {\cal P}}

\def\nAM0{{\nu }_{{\cal A},M_0}}
\def\nAN0{{\nu }_{{\cal A},N_0}}

\def\End{ {\rm End }}

\def\CJ{ {\cal J }}

\def\CP{ {\cal P }}
\def\det{ {\rm det }}

\def\ga{\mathfrak{a}}

\def\gp{\mathfrak{p}}

\def\GL{{\rm GL}}
\def\SL{{\rm SL}}

\def\Hom{{\rm Hom}}

\def\di!{\frac{\der^i}{i!}}
\def\dik!{\frac{\der^k_i}{k!}}

\def\id{{\rm id}}

\def\N{\mathbb{N}}

\def\0{\overline{0}}
\def\1{\overline{1}}

\def\Ln1{\L_{n,\overline{1}}}

\def\a1{a_{\overline{1}}}

\def\S{\Sigma}

\def\sign{{\rm sign}}
\def\vn1{\overrightarrow{n-1}}

\def\Aff{{\rm Aff}}

\def\im{{\rm im}}

\def\mA{\mathbb{A}}

\def\Inn{{\rm Inn}}

\def\mS{\mathbb{S}}
\def\mJ{\mathbb{J}}

\def\clKdim{{\rm cl.Kdim}}

\def\mT{\mathbb{T}}
\def\ind{{\rm ind}}

\def\mE{\mathbb{E}}
\def\mX{\mathbb{X}}

\def\bgp{\overline{\gp}}

\def\K1{{\rm K}_1}
\def\bG{\overline{\G}}
\def\bpsi{\overline{\psi}}
\def\mK{\mathbb{K}}
\def\mP{\mathbb{P}}
\def\mY{\mathbb{Y}}
\def\SAff{{\rm SAff}}
\def\aff{{\rm aff}}
\begin{document}

\author{V. V. \  Bavula 
}

\title{The group of automorphisms of the  algebra
 of one-sided inverses of a polynomial algebra, II }

\maketitle

\begin{abstract}
The algebra $\mS_n$ of one-sided inverses of a polynomial algebra
$P_n$ in $n$ variables is obtained from $P_n$ by adding commuting,
 {\em left} (but not two-sided) inverses of the canonical
 generators of the algebra $P_n$. Ignoring non-Noetherian
 property, the algebra $\mS_n$ belongs to a family of algebras
 like the $n$th Weyl algebra $A_n$ and the polynomial algebra
 $P_{2n}$. Explicit generators are found for the group $G_n$ of
 automorphisms of the algebra $\mS_n$ and for the group $\mS_n^*$
 of units of $\mS_n$ (both groups are huge). An analog of the {\em
Jacobian}  homomorphism  $\CP_n:=\Aut_{K-{\rm alg}}(P_{n} )\ra
K^*$  is introduced for the group $G_n$ (notice that the algebra
$\mS_n$ is noncommutative and neither left nor right Noetherian).
The polynomial Jacobian homomorphism  is {\em unique} since
$\CP_n/[\CP_n, \CP_n]\simeq K^*$. Its analogue is also unique for
$n>2$ but for $n=1,2$ there are exactly two of them. The proof is
based on the following theorem:
$$ G_n/[G_n,G_n]\simeq \begin{cases}
K^*\times K^*& \text{if } n=1,\\
\Z/ 2\Z \times  K^* \times \Z/ 2\Z  &\text{if } n=2, \\
\Z/ 2\Z \times  K^*  &   \text{if } n>2. \\
\end{cases} $$




\end{abstract}


\section{Introduction}
Throughout, ring means an associative ring with $1$; module means
a left module;
 $\N :=\{0, 1, \ldots \}$ is the set of natural numbers; $K$ is a
field of characteristic zero and  $K^*$ is its group of units;
$P_n:= K[x_1, \ldots , x_n]$ is a polynomial algebra over $K$;
$\der_1:=\frac{\der}{\der x_1}, \ldots , \der_n:=\frac{\der}{\der
x_n}$ are the partial derivatives ($K$-linear derivations) of
$P_n$; $\End_K(P_n)$ is the algebra of all $K$-linear maps from
$P_n$ to $P_n$ and $\Aut_K(P_n)$ is its group of units (i.e. the
group of all the invertible linear maps from $P_n$ to $P_n$); the
subalgebra  $A_n:= K \langle x_1, \ldots , x_n , \der_1, \ldots ,
\der_n\rangle$ of $\End_K(P_n)$ is called the $n$'th {\em Weyl}
algebra.

$\noindent $

{\it Definition}, \cite{shrekalg}. The 
{\em algebra} $\mathbb{S}_n$ {\em of one-sided inverses} of $P_n$
is an algebra generated over a field $K$ of characteristic zero by
$2n$ elements $x_1, \ldots , x_n, y_n, \ldots , y_n$ that satisfy
the defining relations:
$$ y_1x_1=1, \ldots , y_nx_n=1, \;\; [x_i, y_j]=[x_i, x_j]= [y_i,y_j]=0
\;\; {\rm for\; all}\; i\neq j,$$ where $[a,b]:= ab-ba$ is  the
algebra commutator of elements $a$ and $b$.

$\noindent $

By the very definition, the algebra $\mS_n$ is obtained from the
polynomial algebra $P_n$ by adding commuting, left (but not
two-sided) inverses of its canonical generators. The algebra
$\mS_1$ is a well-known primitive algebra \cite{Jacobson-StrRing},
p. 35, Example 2.  Over the field
 $\mathbb{C}$ of complex numbers, the completion of the algebra
 $\mS_1$ is the {\em Toeplitz algebra} which is the
 $\mathbb{C}^*$-algebra generated by a unilateral shift on the
 Hilbert space $l^2(\N )$ (note that $y_1=x_1^*$). The Toeplitz
 algebra is the universal $\mathbb{C}^*$-algebra generated by a
 proper isometry.

$\noindent $

{\it Example}, \cite{shrekalg}. Consider a vector space $V=
\bigoplus_{i\in \N}Ke_i$ and two shift operators on $V$, $X:
e_i\mapsto e_{i+1}$ and $Y:e_i\mapsto e_{i-1}$ for all $i\geq 0$
where $e_{-1}:=0$. The subalgebra of $\End_K(V)$ generated by the
operators $X$ and $Y$ is isomorphic to the algebra $\mS_1$
$(X\mapsto x$, $Y\mapsto y)$. By taking the $n$'th tensor power
$V^{\t n }=\bigoplus_{\alpha \in \N^n}Ke_\alpha$ of $V$ we see
that the algebra $\mS_n\simeq \mS_1^{\t n}$ is isomorphic to the
subalgebra of $\End_K(V^{\t n })$ generated by the $2n$ shifts
$X_1, Y_1, \ldots , X_n, Y_n$ that act in different directions.

$\noindent $

The algebra $\mS_n$ is a noncommutative, non-Noetherian algebra
which is not a domain either. Moreover, it contains the algebra of
infinite dimensional matrices. The Gelfand-Kirillov dimension  and
the classical Krull dimension  of the algebra $\mS_n$ is $2n$, but
the global dimension and the weak homological dimension of the
algebra $\mS_n$ is $n$, \cite{shrekalg}.

$\noindent $

{\bf Explicit generators for the group $G_n$}.
 Let $G_n:=\Aut_{K-{\rm
alg}}(\mS_n)$ be the group of automorphisms of the algebra
$\mS_n$, and $\mS_n^*$ be the group of units of the algebra
$\mS_n$. The groups $G_n$ and $\mS_n^*$ are huge, eg both of them
contain the group $\underbrace{\GL_\infty (K)\ltimes\cdots \ltimes
\GL_\infty (K)}_{2^n-1 \;\; {\rm times}}$ which is a small part of
them. A semi-direct product ${}^{semi}\prod_{i=1}^m H_i=H_1\ltimes
H_2\ltimes \cdots \ltimes H_m$ of several groups means that
$H_1\ltimes (H_2\ltimes ( \cdots \ltimes (H_{m-1}\ltimes
H_m)\cdots )$.

\begin{theorem}\label{Int24Apr9}
\begin{enumerate}
\item \cite{shrekaut} $\; G_n=S_n\ltimes \mT^n\ltimes \Inn
(\mS_n)$.\item \cite{shrekaut} $\; G_1\simeq \mT^1\ltimes
\GL_\infty (K)$.
\end{enumerate}
\end{theorem}
In the theorem above,  $S_n=\{ s \in S_n\, | \, s (x_i) =
x_{s(i)}, s(y_i)= y_{s(i)}\}$ is the symmetric group, $\mT^n:=\{
t_\l \, | \, t_\l (x_i) = \l _ix_i, t_\l (y_i) = \l_i^{-1}y_i, \l
=(\l_i)\in K^{*n}\}$ is the $n$-dimensional torus, $\Inn (\mS_n)$
is the group of inner automorphisms of the algebra $\mS_n$ (which
is huge), and $\GL_\infty (K)$ is the group of all the invertible
infinite dimensional matrices of the type $1+M_\infty (K)$ where
the algebra (without 1) of infinite dimensional matrices $M_\infty
(K) :=\varinjlim M_d(K)=\bigcup_{d\geq 1}M_d(K)$ is the injective
limit of matrix algebras.  Theorem \ref{Int24Apr9} is a difficult
one (see the Introduction of \cite{shrekaut} where the structure
and the main ideas of the proof are explained).

The results of the papers \cite{Bav-Jacalg, shrekalg,  shrekaut,
jacaut, K1aut} and of the present paper show that (when ignoring
non-Noetherian property) the algebra $\mS_n$ belongs to a family
of algebras like the $n$'th Weyl algebra $A_n$,  the polynomial
algebra $P_{2n}$, and the Jacobian algebra $\mA_n$ (see
\cite{Bav-Jacalg, jacaut}). Moreover, the algebras $\mS_n$,
$\mA_n$, and $A_n$ are {\em generalized Weyl algebras}.  The
structure of the group $G_1\simeq \mT^1\ltimes \GL_\infty (K)$ is
another confirmation of `similarity' of the algebras $P_2$, $A_1$,
and $\mS_1$. The groups of automorphisms of the polynomial algebra
$P_2$ and the Weyl algebra $A_1$ were found by Jung \cite{jung},
Van der Kulk \cite{kulk}, and Dixmier \cite{Dix} respectively.
These two groups have almost identical structure, they are
`infinite $GL$-groups' in the sense that they are generated by the
torus $\mT^1$ and by the obvious automorphisms: $x\mapsto x+\l
y^i$, $y\mapsto y$; $x\mapsto x$, $y\mapsto y+ \l x^i$, where
$i\in \N$ and $\l \in K$; which are sort of `elementary infinite
dimensional matrices' (i.e. `infinite dimensional transvections`).
The same picture as for the group $G_1$. In prime characteristic,
the group of automorphism of the Weyl algebra $A_1$ was found by
Makar-Limanov \cite{Mak-LimBSMF84} (see also Bavula
\cite{A1rescen}  for a different approach and for further
developments).

A next step in finding explicitly the group $G_n$ and its
generators is done in \cite{K1aut} where explicit generators are
found for the group $G_2$ and it is proved that
\begin{itemize}
\item  (Theorem 2.12, \cite{K1aut}) $ G_2\simeq S_2\ltimes
\mT^2\ltimes \Z\ltimes ((K^*\ltimes E_\infty
(\mS_1))\boxtimes_{\GL_\infty (K)}(K^*\ltimes E_\infty (\mS_1)))$
  {\em where $E_\infty (\mS_1)$ is the subgroup of $\GL_\infty (\mS_1)$
generated by the elementary matrices.}
\end{itemize}
The aim of the present paper is to find explicitly the group $G_n$
(Theorem \ref{B10May9})
 and its generators for $n\geq 2$.

\begin{itemize}
\item  (Theorem \ref{A25May9}) {\em Let $J_s:= \{ 1, \ldots , s\}$
where $s=1, \ldots , n$. The group $G_n= S_n\ltimes \mT^n\ltimes
\Inn (\mS_n)$ is generated by the transpositions $(ij)$ where
$i<j$; the elements $t_{(\l , 1, \ldots , 1)}:x_1\mapsto \l x_1$,
$y_1\mapsto \l^{-1}y_1$, $x_k\mapsto x_k$, $y_k\mapsto y_k$, $k=2,
\ldots , n$; and the inner automorphisms $\o_u$ where $u$ belongs
to the following sets:
\begin{enumerate}
\item $\th_{s,1}(J_s):=
(1+(y_s-1)\prod_{i=1}^{s-1}(1-x_iy_i))\cdot
(1+(x_1-1)\prod_{j=2}^s(1-x_jy_j))$, $s=2, \ldots , n$;
 \item
$1+x_n^tE_{0\alpha}(J_s)$, $1+x_n^tE_{\alpha 0}(J_s)$,
$1+y_n^tE_{0\alpha}(J_s)$, and $1+y_n^tE_{\alpha 0}(J_s)$ where
$t\in \N\backslash \{ 0\}$, $s=1, \ldots , n-1$, and  $\alpha \in
\N^s\backslash \{ 0\}$; and \item $1+(\l -1)E_{00}(J_s)$,
$1+E_{0\alpha} (J_s)$, and $1+E_{\alpha 0}(J_s)$ where $\l \in
K^*$,  $s=1, \ldots , n$, and $\alpha \in \N^s\backslash \{ 0\}$.
\end{enumerate}
where $E_{00}(J_s):= \prod_{i=1}^s(1-x_iy_i)$,  $E_{0\alpha} (J_s)
:=\prod_{i=1}^s(y_i^{\alpha_i}-x_iy_i^{\alpha_i+1})$,  and
$E_{\alpha 0} (J_s)
:=\prod_{i=1}^s(x_i^{\alpha_i}-x_i^{\alpha_i+1}y_i)$. }
\end{itemize}

{\bf The structure and main ideas of finding the generators for
the groups $G_n$ and $\mS_n^*$}.  A first step is the following
theorem which is a consequence of an analogous result for the
Jacobian algebra $\mA_n$ (Theorem 4.4, \cite{Bav-Jacalg}) modulo a
discrete subgroup. Then the theorem follows from the inclusion
$\mS^*\subseteq \mA_n^*$.

\begin{theorem}\label{aInt24Apr9}
\cite{jacaut}
\begin{enumerate}
\item $\mS_n^* = K^*\times (1+\ga_n)^*$ where the ideal $\ga_n$ of
the algebra $\mS_n$ is the sum of all the height one prime ideals
of the algebra $\mS_n$.  \item  The centre of the group $\mS_n^*$
is $K^*$, and the centre of the group $(1+\ga_n)^*$ is $\{ 1 \}$.
\item The map $(1+\ga_n)^*\ra \Inn (\mS_n)$, $u\mapsto \o_u$, is a
group  isomorphism.
\end{enumerate}
\end{theorem}
This theorem reduces the problem of finding the group $G_n$ to the
problem of finding the group of units $(1+\ga_n)^*$. To save on
notation, often we identify the groups $(1+\ga_n)^*$ and $\Inn
(\mS_n)$ via $u\mapsto \o_u$.

The polynomial algebra $P_n$ is a faithful $\mS_n$-module (see
Example above), hence $\mS_n\subset \End_K(P_n)$. The ideals of
the algebra $\mS_n$ commute ($IJ=JI$), \cite{shrekalg}. There are
precisely $n$ height one prime ideals of the algebra $\mS_n$, say
$\gp_1, \ldots , \gp_n$. They are found explicitly in
\cite{shrekalg}, and they form a single $G_n$-orbit. In
particular, the ideals $\ga_{n,s}:=\sum_{i_1<\cdots
<i_s}\gp_{i_1}\cdots \gp_{i_s}$, $s=1, \ldots , n$, are
$G_n$-invariant ideals of the algebra $\mS_n$. The group
$(1+\ga_n)^*$ has the strictly descending chain of
$G_n$-invariant (hence normal) subgroups
$$ (1+\ga_n)^*=(1+\ga_{n,1})^*\supset \cdots \supset
(1+\ga_{n,s})^*\supset \cdots \supset (1+\ga_{n,n-1})^*\supset
(1+\ga_{n,n})^*.$$ Briefly, to prove results on the group
$(1+\ga_n)^*$ we first prove similar results for the subgroups
$(1+\ga_{n,s})^*$, $s=1, \ldots , n-1$, using a double induction
on $(n,s)$ starting with $(n,n-1)$ in the second part of the
induction (the induction on $s$ is a downward induction, the group
$(1+\ga_{n,n})^*$ is isomorphic to $\GL_\infty (K)$ and contains
no essential information about the overgroups, that is why we have
to start with $(n,n-1)$). The initial case $(n,n-1)$ is the most
difficult one, we spend entire Section \ref{GTNN1} to treat it.

Difficulty in finding the group $(1+\ga_n)^*$ stems from two
facts:  (i) $\mS_n^*\subsetneqq \mS_n\cap \Aut_K(P_n)$, i.e. there
are non-units of the algebra $\mS_n$ that are invertible linear
maps in $P_n$; and (ii) some units of the algebra $\mS_n$ are
product of  {\em non-units}. To tackle  the second problem the
so-called  {\em current groups} $\Theta_{n,s}$, $s=1, \ldots,
n-1$, are introduced. These are finitely generated subgroups of
$(1+\ga_n)^*$ generated by explicit generators, and each of the
generators is a product of two non-units of $\mS_n^*$ (they are
even non-units of $\End_K(P_n)$). The current groups turn out to
be the most important subgroups of $(1+\ga_n)^*$, they control the
most difficult parts of the structure of the group $(1+\ga_n)^*$.

In dealing with the case $(n,n-1)$, we use the Fredholm linear
maps/operators and their indices. This technique is not available
in other cases, i.e. when $(n,s)\neq (n,n-1)$, but the point is
that other cases can be reduced to  the initial one but over a
larger coefficient {\em ring} (not a field). The indices of
operators are used to construct several group homomorphisms. The
most difficult part of Section \ref{GTNN1} is to prove that the
homomorphisms are well defined maps (as their constructions are
based on highly non-unique decompositions). As a result the group
$(1+\ga_n)^*$ is found explicitly.

\begin{itemize}
\item (Theorem \ref{B10May9}) $(1+\ga_n)^* = \Theta_{n,1}'
\mE_{n,1} \Theta_{n,2}' \mE_{n,2}\cdots \Theta_{n,n-1}'
\mE_{n,n-1} $,
\end{itemize}
where the sets $\Theta_{n,s}'\subseteq \Theta_{n,s}$ and the
groups $\mE_{n,s}$   are given explicitly, (\ref{Tpns}) and
(\ref{mEns}). As a consequence we have explicit generators for the
group $(1+\ga_n)^*$.

\begin{itemize}
\item (Theorem \ref{25May9}) {\em The group $(1+\ga_n)^*$ is
generated by the following elements:
\begin{enumerate}
\item $\th_{\max (J), j}(J):=(1+(y_{\max (J)}-1)\prod_{i\in
J\backslash \max (J)}(1-x_iy_i))\cdot (1+(x_j-1)\prod_{k\in
J\backslash j}(1-x_ky_k))$ where $J$ runs through all the subsets
of $\{ 1,\ldots , n\}$ that contain at least two elements,  $j\in
J\backslash \max (J)$; \item $1+x_i^tE_{0\alpha}(I)$,
$1+x_i^tE_{\alpha 0}(I)$, $1+y_i^tE_{0\alpha}(I)$, and
$1+y_i^tE_{\alpha 0}(I)$ where  $I$ runs through all the subsets
of $\{ 1,\ldots , n\}$ such that  $|I|=1, \ldots , n-1$, $t\in
\N\backslash \{ 0\}$, $i\not\in I$, $\alpha \in \N^I\backslash \{
0\}$; and \item $1+(\l -1)E_{00}(I)$, $1+E_{0\alpha} (I)$, and
$1+E_{\alpha 0}(I)$ where $\l \in K^*$, $I\neq \emptyset$, and
$\alpha \in \N^I\backslash \{ 0\}$.
\end{enumerate} }
\end{itemize}
Then it is easy to obtain explicit generators for the group $G_n$
(Theorem \ref{A25May9}).

$\noindent $

{\bf An analogue of the polynomial Jacobian homomorphism for the
group $G_n$}. For the polynomial algebra $P_n$ there is an
important group homomorphism, 
\begin{equation}\label{polJac}
 \CJ_n : \CP_n:=\Aut_{K-{\rm alg}}(P_{n})\ra K^*, \;\; \s \mapsto
\det (\frac{\der \s (x_i)}{\der x_j}),
\end{equation}
the so-called {\em Jacobian} homomorphism. Note that the Jacobian
homomorphism is a determinant. Each automorphism $\s \in \CP_n$ is
a unique product $\s_{aff}\xi$ of an affine automorphism
$\s_{aff}\in \Aff_n$ and an element $\xi$ of the {\em Jacobian}
group $\S_n$ (see Section \ref{COMGJAC} for details), and the
Jacobian of $\s$ is uniquely determined by its affine part, i.e.
$\CJ (\s ) = \CJ (\s_{aff})$. This property {\em uniquely}
characterizes the Jacobian homomorphism since $$\CP_n/[\CP_n,
\CP_n]\simeq \Aff_n/\Aff_n\cap [\CP_n, \CP_n]\simeq K^*.$$ So,
there are two different ways of defining the Jacobian
homomorphism: by the explicit formula (\ref{polJac}) or as a group
homomorphism from $\CP_n/ [ \CP_n, \CP_n]$ to $K^*$ that is
defined naturally (i.e. as the determinant)  on the affine
subgroup $\Aff_n$ of $\CP_n$.

The group $G_n=S_n\ltimes \mT^n\ltimes \Inn (\mS_n)$ has a similar
structure as the group $\CP_n$ where $\aff_n:=S_n\ltimes \mT^n$ is
an affine part of $G_n$ and the group $\Inn (\mS_n)$ of inner
automorphisms plays a role of the Jacobian group. In Section
\ref{COMGJAC}, we introduce an analogue $\mJ_n : G_n\ra K^*$ of
the Jacobian homomorphism  using the second definition of the
polynomial Jacobian map as a guiding principle: the map $\mJ_n$ is
a homomorphism $\mJ_n : G_n/ [G_n,G_n]\ra K^*$ such that on the
affine group $\aff_n$ it is defined exactly in the same way as in
the polynomial case. For $n>2$, the homomorphism $\mJ_n$ is {\em
unique} (Theorem \ref{15Jun9}) by the same reason as in the
polynomial case, namely, (Corollary \ref{a12Jun9}.(3))
$$G_n/[G_n,G_n]\simeq \aff_n/[\aff_n, \aff_n]\simeq \Z / 2\Z \times  K^*.$$
But for $n=1,2$, the homomorphism $\mJ_n$ is not unique, there are
exactly two of them since  (Corollary \ref{a12Jun9}.(3))
$$G_n/[G_n,G_n]\simeq \aff_n/[\aff_n, \aff_n] \times \begin{cases}
K^*& \text{if } n=1,\\
\Z / 2\Z & \text{if }n=2.\\
\end{cases}
$$
More informally, for $n=1$ the appearance of the second, exotic
$\mJ_1^{ex}$ has connection with existence of the determinant
(homomorphism) on the group $\GL_\infty (K)$, but for $n=2$ the
exotic $\mJ_2^{ex}$ is explained by the fact that the current
group $\Theta_2$ does not belong to the commutant $[G_2,G_2]$. For
$n=1$, $\mJ_1$ and $\mJ_1^{ex}$ are algebraically independent
characters of the group $G_1$, but, for $n=2$, $\mJ_2^2=
(\mJ_2^{ex})^2$ (Theorem \ref{15Jun9}).

 The proofs are  based on finding explicitly
the commutant $[G_n,G_n]$ of the group $G_n$ (Theorem
\ref{6Jun9}.(1)) and proving that
\begin{itemize}
\item (Theorem \ref{6Jun9}.(2))$$ G_n/[G_n,G_n]\simeq
\begin{cases}
K^*\times K^*& \text{if } n=1,\\
\Z/ 2\Z \times  K^* \times \Z/ 2\Z  &\text{if } n=2, \\
\Z/ 2\Z \times  K^* &   \text{if } n>2. \\
\end{cases}$$
\end{itemize}
The most surprising thing is that despite the fact that the
algebra $\mS_n$ is noncommutative, non-Noetherian, of
Gelfand-Kirillov dimension $2n$ (not $n$), not a domain, the
unique  homomorphism $\mJ_n$ `coincides' with the polynomial
Jacobian homomorphism $\CJ_n$ for the polynomial algebra $P_n$
(not $P_{2n}$): for $\s \in G_n$, 
\begin{equation}\label{SnJac}
\mJ_n (\s ) = \det (\frac{\der \overline{\s } (x_i)}{\der x_j}),
\end{equation}
i.e. the homomorphism $\mJ_n$ is the composition of two
homomorphisms:
$$ G_n\ra \Aut_{K-{\rm alg}}(\mS_n/\ga_n\simeq L_n), \;\; \s
\mapsto \overline{\s}: a+\ga_n\mapsto \s (a) +\ga_n, $$ where
$L_n:=K[x_1,x_1^{-1}, \ldots , x_n, x_n^{-1}]$ is the Laurent
polynomial algebra (see (\ref{SnSn})) and
$$ \CJ_n :  \Aut_{K-{\rm alg}}( L_n)\ra L_n^*, \;\; \tau \mapsto \det (\frac{\der \tau  (x_i)}{\der
x_j}).$$ {\bf Proof of (\ref{SnJac})}. Since $G_n=\aff_n\ltimes
\Inn (\mS_n)$ and the factor algebra $\mS_n / \ga_n \simeq L_n$ is
{\em commutative}, the homomorphism  (\ref{SnJac}) acts trivially
on $\Inn (\mS_n)$ (i.e. $\mJ_n ( \Inn (\mS_n))=1$), but on
$\aff_n$ the map (\ref{SnJac}) acts exactly as in the polynomial
case: for each element $s\cdot t_\l \in \aff_n$ where $s\in S_n$
and $t_\l \in \mT^n$, $\mJ_n(s\cdot t_\l )=\sign
(s)\cdot\prod_{i=1}^n\l_i$ where $\sign (s)\in \{ \pm 1\}$ is the
sign/parity of the permutation $s$. $\Box$

$\noindent $

So, for each  element $\s = s\cdot t_\l \cdot \o_u\in G_n =
S_n\ltimes \mT^n\ltimes \Inn (\mS_n)$ where $s\in S_n$, $t_\l \in
\mT^n$, and $\o_u\in \Inn (\mS_n)$, 
\begin{equation}\label{SnJac1}
\mJ_n (\s ) =\sign (s)\cdot\prod_{i=1}^n\l_i.
\end{equation}
One may have noticed that  $\mJ_n (\s )$ depends only on $\s
(x_1),\ldots , \s (x_n)$ and the set $\{ x_1, \ldots , x_n\}$ is
not a generating set for the algebra $\mS_n$. It is a trivial
observation that an algebra endomorphism is  uniquely determined
by its action on a set of algebra generators but, for the algebra
$\mS_n$, an algebra endomorphism is uniquely determined by its
action on either of the sets $\{ x_1, \ldots , x_n\}$  or $\{ y_1,
\ldots , y_n\}$ (which are not algebra generating sets).

\begin{theorem}\label{6Feb9}
{\rm (Rigidity of the group
$G_n$.  Theorem 3.7, \cite{shrekaut})}  Let $\s , \tau \in G_n$.
Then the following statements are equivalent.
\begin{enumerate}
 \item $\s = \tau$. \item $\s (x_1) = \tau (x_1), \ldots , \s (x_n) = \tau
 (x_n)$.
 \item$\s (y_1) = \tau (y_1), \ldots , \s (y_n) = \tau (y_n)$.
\end{enumerate}
\end{theorem}


\section{The group $(1+\ga_n)^*$ and its subgroups} \label{KTG1S}

In this section, we collect some results without proofs on the
algebras $\mS_n$ from \cite{shrekalg} that will be used in this
paper, their proofs can be found in \cite{shrekalg}.
 Several important
subgroups of the group $(1+\ga_n)^*$ are introduced. The most
interesting of these are the current subgroups $\Theta_{n,s}$,
$s=1, \ldots , n-1$. They encapsulate the most difficult parts of
the groups $\mS_n^*$ and $G_n$.  This section sets a scene for
proving the main results of the paper.

$\noindent $

{\bf The algebra of one-sided inverses of a polynomial algebra}.
Clearly, $\mathbb{S}_n=\mS_1(1)\t \cdots \t\mS_1(n)\simeq
\mathbb{S}_1^{\t n}$ where $\mS_1(i):=K\langle x_i,y_i \, | \,
y_ix_i=1\rangle \simeq \mS_1$ and $$\mS_n=\bigoplus_{\alpha ,
\beta \in \N^n} Kx^\alpha y^\beta$$ where $x^\alpha :=
x_1^{\alpha_1} \cdots x_n^{\alpha_n}$, $\alpha = (\alpha_1, \ldots
, \alpha_n)$, $y^\beta := y_1^{\beta_1} \cdots y_n^{\beta_n}$,
$\beta = (\beta_1,\ldots , \beta_n)$. In particular, the algebra
$\mS_n$ contains two polynomial subalgebras $P_n$ and $Q_n:=K[y_1,
\ldots , y_n]$ and is equal,  as a vector space,  to their tensor
product $P_n\t Q_n$. Note that also the Weyl algebra $A_n$ is a
tensor product (as a vector space) $P_n\t K[\der_1, \ldots ,
\der_n]$ of its two polynomial subalgebras.

When $n=1$, we usually drop the subscript `1' if this does not
lead to confusion.  So, $\mS_1= K\langle x,y\, | \,
yx=1\rangle=\bigoplus_{i,j\geq 0}Kx^iy^j$. For each natural number
$d\geq 1$, let $M_d(K):=\bigoplus_{i,j=0}^{d-1}KE_{ij}$ be the
algebra of $d$-dimensional matrices where $\{ E_{ij}\}$ are the
matrix units, and
$$M_\infty (K) :=
\varinjlim M_d(K)=\bigoplus_{i,j\in \N}KE_{ij}$$ be the algebra
(without 1) of infinite dimensional matrices. The algebra $\mS_1$
contains the ideal $F:=\bigoplus_{i,j\in \N}KE_{ij}$, where
\begin{equation}\label{Eijc}
E_{ij}:= x^iy^j-x^{i+1}y^{j+1}, \;\; i,j\geq 0.
\end{equation}
For all natural numbers $i$, $j$, $k$, and $l$,
$E_{ij}E_{kl}=\d_{jk}E_{il}$ where $\d_{jk}$ is the Kronecker
delta function.  The ideal $F$ is an algebra (without 1)
isomorphic to the algebra $M_\infty (K)$ via $E_{ij}\mapsto
E_{ij}$.  For all $i,j\geq 0$, 
\begin{equation}\label{xyEij}
xE_{ij}=E_{i+1, j}, \;\; yE_{ij} = E_{i-1, j}\;\;\; (E_{-1,j}:=0),
\end{equation}
\begin{equation}\label{xyEij1}
E_{ij}x=E_{i, j-1}, \;\; E_{ij}y = E_{i, j+1} \;\;\;
(E_{i,-1}:=0).
\end{equation}
The algebra
\begin{equation}\label{mS1d}
\mS_1= K\oplus xK[x]\oplus yK[y]\oplus F
\end{equation}
is the direct sum of vector spaces. Then 
\begin{equation}\label{mS1d1}
\mS_1/F\simeq K[x,x^{-1}]=:L_1, \;\; x\mapsto x, \;\; y \mapsto
x^{-1},
\end{equation}
since $yx=1$, $xy=1-E_{00}$ and $E_{00}\in F$.

$\noindent $

The algebra $\mS_n = \bigotimes_{i=1}^n \mS_1(i)$ contains the
ideal
$$F_n:= F^{\t n }=\bigoplus_{\alpha , \beta \in
\N^n}KE_{\alpha \beta}, \;\; {\rm where}\;\; E_{\alpha
\beta}:=\prod_{i=1}^n E_{\alpha_i \beta_i}(i), \;\; E_{\alpha_i
\beta_i}(i):=x_i^{\alpha_i}y_i^{\beta_i}-x_i^{\alpha_i+1}y_i^{\beta_i+1}.$$
Note that $E_{\alpha \beta}E_{\g \rho}=\d_{\beta \g }E_{\alpha
\rho}$ for all elements $\alpha, \beta , \g , \rho \in \N^n$ where
$\d_{\beta
 \g }$ is the Kronecker delta function.

{\bf The involution $\eta$ on $\mS_n$}. The algebra $\mS_n$ admits
the {\em involution}
$$ \eta : \mS_n\ra \mS_n, \;\; x_i\mapsto y_i, \;\; y_i\mapsto
x_i, \;\; i=1, \ldots , n,$$ i.e. it is a $K$-algebra
anti-isomorphism ($\eta (ab) = \eta (b) \eta (a)$ for all $a,b\in
\mS_n$) such that $\eta^2 = \id_{\mS_n}$, the identity map on
$\mS_n$. So, the algebra $\mS_n$ is {\em self-dual} (i.e. it is
isomorphic to its opposite algebra, $\eta : \mS_n\simeq
\mS_n^{op}$). The involution $\eta$ acts on the `matrix' ring
$F_n$ as the transposition,  
\begin{equation}\label{eEij1}
\eta (E_{\alpha \beta} )=E_{\beta \alpha}.
\end{equation}
 The canonical generators $x_i$,
$y_j$ $(1\leq i,j\leq n)$ determine the ascending filtration $\{
\mS_{n, \leq i}\}_{i\in \N}$ on the algebra $\mS_n$ in the obvious
way (i.e. by the total degree of the generators): $\mS_{n, \leq
i}:= \bigoplus_{|\alpha |+|\beta |\leq i} Kx^\alpha y^\beta$ where
$|\alpha | \; = \alpha_1+\cdots + \alpha_n$ ($\mS_{n, \leq
i}\mS_{n, \leq j}\subseteq \mS_{n, \leq i+j}$ for all $i,j\geq
0$). Then $\dim (\mS_{n,\leq i})={i+2n \choose 2n}$ for $i\geq 0$,
and so the Gelfand-Kirillov dimension $\GK (\mS_n )$ of the
algebra $\mS_n$ is equal to $2n$. It is not difficult to show
 that the algebra $\mS_n$ is neither left nor
right Noetherian. Moreover, it contains infinite direct sums of
left and right ideals (see \cite{shrekalg}).

\begin{itemize}
 \item   {\em The
 algebra $\mS_n$ is central, prime, and catenary. Every nonzero
 ideal  of $\mS_n$ is an essential left and right submodule of}
 $\mS_n$.
 \item  {\em The ideals of
 $\mS_n$ commute ($IJ=JI$);  and the set of ideals of $\mS_n$ satisfy the a.c.c..}
 \item  {\em The classical Krull dimension $\clKdim (\mS_n)$ of $\mS_n$ is $2n$.}
  \item  {\em Let $I$
be an ideal of $\mS_n$. Then the factor algebra $\mS_n / I$ is
left (or right) Noetherian iff the ideal $I$ contains all the
height one primes of $\mS_n$.}
\end{itemize}

{\bf The set of height 1 primes of $\mS_n$}.   Consider the ideals
of the algebra $\mS_n$:
$$\gp_1:=F\t \mS_{n-1},\; \gp_2:= \mS_1\t F\t \mS_{n-2}, \ldots ,
 \gp_n:= \mS_{n-1} \t F.$$ Then $\mS_n/\gp_i\simeq
\mS_{n-1}\t (\mS_1/F) \simeq  \mS_{n-1}\t K[x_i, x_i^{-1}]$ and
$\bigcap_{i=1}^n \gp_i = \prod_{i=1}^n \gp_i =F^{\t n }$. Clearly,
$\gp_i\not\subseteq \gp_j$ for all $i\neq j$.

\begin{itemize}
 \item
{\em The set $\CH_1$ of height one prime ideals of the algebra
$\mS_n$ is} $\{ \gp_1, \ldots , \gp_n\}$.
\end{itemize}

 Let
$\ga_n:= \gp_1+\cdots +\gp_n$. Then the factor algebra
\begin{equation}\label{SnSn}
\mS_n/ \ga_n\simeq (\mS_1/F)^{\t n } \simeq \bigotimes_{i=1}^n
K[x_i, x_i^{-1}]= K[x_1, x_1^{-1}, \ldots , x_n, x_n^{-1}]=:L_n
\end{equation}
is a skew Laurent polynomial algebra in $n$ variables,  and so
$\ga_n$ is a prime ideal of height and co-height $n$ of the
algebra $\mS_n$.

\begin{proposition}\label{a19Dec8}
\cite{shrekalg} The polynomial algebra $P_n$
 is the only faithful, simple $\mS_n$-module.
\end{proposition}

In more detail, ${}_{\mS_n}P_n\simeq \mS_n / (\sum_{i=0}^n \mS_n
y_i) =\bigoplus_{\alpha \in \N^n} Kx^\alpha \overline{1}$,
$\overline{1}:= 1+\sum_{i=1}^n \mS_ny_i$; and the action of the
canonical generators of the algebra $\mS_n$ on the polynomial
algebra $P_n$ is given by the rule:
$$ x_i*x^\alpha = x^{\alpha + e_i}, \;\; y_i*x^\alpha = \begin{cases}
x^{\alpha - e_i}& \text{if } \; \alpha_i>0,\\
0& \text{if }\; \alpha_i=0,\\
\end{cases}  \;\; {\rm and }\;\; E_{\beta \g}*x^\alpha = \d_{\g
\alpha} x^\beta,
$$
where $e_1:= (1,0,\ldots , 0),  \ldots , e_n:=(0, \ldots , 0,1)$
is the canonical basis for the free $\Z$-module
$\Z^n=\bigoplus_{i=1}^n \Z e_i$.  We identify the algebra $\mS_n$
with its image in the algebra $\End_K(P_n)$ of all the $K$-linear
maps from the vector space $P_n$ to itself, i.e. $\mS_n \subset
\End_K(P_n)$.

For each non-empty subset $I$ of the set $\{ 1,\ldots , n\}$, let
$\mS_I:=\bigotimes_{i\in I}\mS_1(i)\simeq \mS_{|I|}$  where $|I|$
is the number of elements in the set $I$, $F_I:=\bigotimes_{i\in
I}F(i)\simeq M_\infty (K)$, $\ga_I$ be the ideal of the algebra
$\mS_I$ generated by the vector space $\bigoplus_{i\in I}F(i)$,
i.e. $\ga_I:=\sum_{i\in I}F(i)\t \mS_{I\backslash i}$. The factor
algebra $L_I:=\mS_I/\ga_I\simeq K[x_i,x_i^{-1}]_{i\in I}$ is a
Laurent polynomial algebra. For elements $\alpha =
(\alpha_i)_{i\in I}, \beta = (\beta_i)_{i\in I}\in \N^I$, let $
E_{\alpha \beta}:=\prod_{i\in I}E_{\alpha_i\beta_i}(I)$. Then
$E_{\alpha\beta}E_{\xi\rho}= \d_{\beta \xi}E_{\alpha\rho}(I)$ for
all $\alpha , \beta , \xi , \rho \in \N^I$.

$\noindent $

{\bf The $G_n$-invariant normal subgroups $(1+\ga_{n,s})^*$  of
$(1+\ga_n)^*$.} We will use often the following two obvious
lemmas.

\begin{lemma}\label{aa13Dec8}
\cite{shrekaut} Let $R$ be a ring and $I_1, \ldots , I_n$ be
ideals of the ring $R$ such that $I_iI_j=0$ for all $i\neq j$. Let
$a= 1+a_1+\cdots +a_n\in R$ where $a_1\in I_1, \ldots , a_n\in
I_n$. The element $a$ is a unit of the ring $R$ iff all the
elements $1+a_i$ are units; and, in this case,  $a^{-1}=
(1+a_i)^{-1} (1+a_2)^{-1}\cdots (1+a_n)^{-1}$.
\end{lemma}

Let $R$ be a ring, $R^*$ be its group of units, $I$ be an ideal of
$R$ such that $I\neq R$, and let  $ (1+I)^*$ be the group of units
of the multiplicative monoid $1+I$.

\begin{lemma}\label{a14Dec8}
\cite{shrekaut}  Let $R$ and  $I$ be as above. Then
\begin{enumerate}
\item $R^*\cap (1+I)= (1+I)^*$.\item $(1+I)^*$ is a normal
subgroup of $R^*$.
\end{enumerate}
\end{lemma}

For each subset $I$ of the set $\{ 1, \ldots , n\}$, let $\gp_I:=
\bigcap_{i\in I}\gp_i$, and $\gp_\emptyset :=\mS_n$. Each $\gp_I$
is an ideal of the algebra $\mS_n$ and $\gp_I=\prod_{i\in
I}\gp_i$.  The complement to the subset $I$ is denoted by $CI$.
For an one-element subset $\{ i\}$, we write $Ci$ rather than $C\{
i\}$. In particular, $\gp_{Ci}:= \gp_{C\{ i\} }=\bigcap_{j\neq
i}\gp_j$.

 For each number $s=1,\ldots , n$, let
$\ga_{n,s}:=\sum_{|I|=s}\gp_I$. By the very definition, the ideals
$\ga_{n,s}$ are $G_n$-invariant ideals (since the set $\CH_1$ of
all the height one prime ideals of the algebra $\mS_n$ is $\{
\gp_1, \ldots , \gp_n\}$, \cite{shrekaut}, and $\CH_1$ is a
$G_n$-orbit). We have the strictly descending chain  of
$G_n$-invariant  ideals of the algebra $\mS_n$:
$$ \ga_n=\ga_{n,1}\supset \ga_{n,2}\supset \cdots \supset
\ga_{n,s}\supset \cdots \supset \ga_{n,n}=F_n\supset
\ga_{n,n+1}:=0.$$ These are also ideals of the subalgebra
$K+\ga_n$ of $\mS_n$. Each set $\ga_{n,s}$ is an ideal of the
algebra $K+\ga_{n,t}$ for all $t\leq s$, and the group of units of
the algebra $K+\ga_{n,s}$ is the direct product of its two
subgroups (Lemma \ref{a14Dec8}.(1)),
$$ (K+\ga_{n,s})^* = K^*\times (1+\ga_{n,s})^*, \;\; s=1, \ldots ,
n.$$ The groups $(K+\ga_{n,s})^*$ and $(1+\ga_{n,s})^*$ are
$G_n$-invariant. For each number $s=1, \ldots , n$, the factor
algebra
$$(K+\ga_{n,s})/\ga_{n,s+1}=K\bigoplus
\bigoplus_{|I|=s}\bgp_I$$ contains the idempotent ideals
$\bgp_I:=(\gp_I+\ga_{n,s+1})/\ga_{n,s+1}$ such that $\bgp_I
\bgp_J=0$ for all $I\neq J$ such that $|I|=|J|=s$.

Recall that for a Laurent polynomial algebra $L=K[x,x^{-1}]$,
${\rm K}_1(L)\simeq L^*$, \cite{Swan}, \cite{Bass-book-K-theory},
\cite{Milnor-book-K-theory}, 
\begin{equation}\label{GLLU}
\GL_\infty (L)=U(L)\ltimes E_\infty (L)
\end{equation}
where $E_\infty (L)$ is the subgroup of $\GL_\infty (L)$ generated
by all the {\em elementary matrices} $\{ 1+aE_{ij}\, | \, a\in L,
i,j\in \N, i\neq j\}$, and $U(L):=\{\mu (u):= uE_{00}+1-E_{00}\, |
\, u\in L^*\}\simeq L^*$, $\mu (u)\lra u$. The group $E_\infty
(L)$ is a normal subgroup of $\GL_\infty (L)$, this is true for an
arbitrary coefficient ring.

By Lemma \ref{aa13Dec8} and (\ref{GLLU}), the group of units of
the algebra $(K+\ga_{n,s})/\ga_{n,s+1}=:K+\ga_{n,s}/\ga_{n,s+1}$
is the direct product of groups,
$$ (K+\ga_{n,s}/\ga_{n,s+1})^*=K^*\times
\prod_{|I|=s}(1+\bgp_I)^*\simeq K^*\times \prod_{|I|=s}\GL_\infty
(L_{CI})\simeq K^*\times \prod_{|I|=s}U(L_{CI})\ltimes E_\infty
(L_{CI})$$ since $(1+\bgp_I)^*\simeq (1+M_\infty
(L_{CI}))^*=\GL_\infty (L_{CI})$ where
$L_{CI}:=\mS_{CI}/\ga_{CI}=\bigotimes_{i\in CI}K[x_i, x_i^{-1}]$
is the Laurent polynomial algebra.
 In more detail, for each non-empty subset $I$ of $\{ 1, \ldots ,
 n\}$, let $\Z^I:=\bigoplus_{i\in I}\Z e_i$, it is a subgroup of
 $\Z^n=\bigoplus_{i=1}^n\Z e_i$. Similarly, $\N^I:= \bigoplus_{i\in I}\N e_i$. By (\ref{GLLU}),
\begin{equation}\label{apns}
(1+\bgp_I)^*=U(L_{CI})\ltimes E_\infty (L_{CI})=(U_I(K)\times
\mX_{CI})\ltimes E_\infty (L_{CI})
\end{equation}
where
\begin{eqnarray*}
 U(L_{CI})&:=&\{\mu_I(u):= uE_{00}(I)+1-E_{00}(I)\, | \, u\in
L^*_{CI}\}\simeq L_{CI}^*,\; \mu_I(u)\lra u, \\
L^*_{CI}&=&\{\l x^\alpha \, | \, \l \in K^*, \alpha \in \Z^{CI}\},\\
U_I(K)&:=&\{ \mu_I(\l ):=\l E_{00}(I)+1-E_{00}(I)\, |\, \l \in
K^*\} \simeq
K^*,\;  \mu_I(\l )\lra \l, \\
\mX_{CI}&:=&\{ \mu_I(x^\alpha ):=x^\alpha E_{00}(I)+1-E_{00}(I)\,
| \, \alpha \in
\Z^{CI}\}\simeq \Z^{CI}\simeq \Z^{n-s},\;  \mu_I(x^\alpha )\lra \alpha , \\
E_\infty (L_{CI})&:=&\langle 1+aE_{\alpha\beta}(I)\, | \, a\in
L_{CI}, \alpha,\beta \in \N^I, \alpha\neq \beta \rangle.
\end{eqnarray*}

The algebra  epimorphism  $\psi_{n,s}: K+\ga_{n,s}\ra
(K+\ga_{n,s})/\ga_{n,s+1}$, $ a\mapsto a+\ga_{n,s+1}$, yields the
group homomorphism of their groups of units $(K+\ga_{n,s})^*\ra
(K+\ga_{n,s}/\ga_{n,s+1})^*$ and the kernel of which is
$(1+\ga_{n,s+1})^*$. As a result we have the exact sequence of
group homomorphisms:
$$
\xymatrix{1\ar[r] & (1+\ga_{n,s+1})^*\ar[r]\ar[d]^{=}  & (K+\ga_{n,s})^* \ar[r]\ar[d]^{=} & (K+\ga_{n,s}/\ga_{n,s+1})^* \ar[d]^{=} & \\
 & (1+\ga_{n,s+1})^*  & K^*\times (1+\ga_{n,s})^*  & K^*\times \prod_{|I|=s}(1+\bgp_I)^*  & }$$
which yields the exact sequence of group homomorphisms:
\begin{equation}\label{apns1}
1\ra (1+\ga_{n,s+1})^*\ra (1+\ga_{n,s})^* \stackrel{\psi_{n,s}}\ra
\prod_{|I|=s}(1+\bgp_I)^* \simeq \prod_{|I|=s}\GL_\infty
(L_{CI})\ra \CZ_{n,s}\ra 1.
\end{equation}
Intuitively, the group $\CZ_{n,s}$ represents `relations' that
determine the image $\im (\psi_{n,s})$ as the subgroup of
 $\prod_{|I|=s}(1+\bgp_I)^*$. We will see later that the group
$\CZ_{n,s}$ is a free abelian group of rank ${n\choose s+1}$
(Corollary \ref{a30May9}). So, the image of the map $\psi_{n,s}$
is large.  Note that $\ga_{n,s+1}$ and $\gp_I$ (where $|I|=s)$ are
ideals of the algebra $K+\ga_{n,s}$. By Lemma \ref{a14Dec8}, the
groups $(1+\ga_{n,s+1})^*$ and $(1+\gp_I)^*$ (where $|I|=s$) are
normal subgroups of $(1+\ga_{n,s})^*$. Then the  subgroup
$\Upsilon_{n,s}$ of $(1+\ga_{n,s})^*$ generated by these normal
subgroups is a normal subgroup of $(1+\ga_{n,s})^*$. As a subset
of $(1+\ga_{n,s})^*$, the group $\Upsilon_{n,s}$ is equal to the
product of the groups $(1+\ga_{n,s+1})^*$, $(1+\gp_I)^*$, $|I|=s$,
in {\em arbitrary} order (by their normality), i.e.
\begin{equation}\label{UPns}
\Upsilon_{n,s}=\prod_{|I|=s}(1+\gp_I)^*\cdot (1+\ga_{n,s+1})^*.
\end{equation}
By Theorem \ref{Int24Apr9} and Theorem \ref{aInt24Apr9}, the group
$\Upsilon_{n,s}$ is a $G_n$-invariant (hence, normal) subgroup of
$\mS_n^*$.  We will see that the factor group $(1+\ga_{n,s})^*/
\Upsilon_{n,s}$ is a free abelian group of rank ${n\choose s+1}s$,
see (\ref{Dns1}).

By (\ref{apns}), the direct product of groups
$\prod_{|I|=s}(1+\bgp_I)^*=\mX_{n,s}\ltimes \bG_{n,s}$ is the
semi-direct product of its two subgroups 
\begin{equation}\label{apns2}
\mX_{n,s}:=\prod_{|I|=s}\mX_{CI}\simeq \Z^{{n\choose s}(n-s)}\;\;
{\rm and}\;\;  \bG_{n,s}:=\prod_{|I|=s}U_I(K)\ltimes
E_\infty(L_{CI}).
\end{equation}
For each subset $I$ of $\{ 1,\ldots , n\}$ such that $|I|=s$,
$U_I(K)\ltimes E_\infty (\mS_{CI})$ is the subgroup of
$(1+\gp_I)^*$ where 
\begin{equation}\label{apns3}
U_I(K):=\{ \mu_I(\l )\, | \, \l \in K^*\} \simeq K^*, \; E_\infty
(\mS_{CI}):=\langle 1+aE_{\alpha\beta}(I)\, | \, a\in \mS_{CI},
\alpha \neq  \beta \in \N^I\rangle ,
\end{equation}
where $\mu_I(\l ):=\l E_{00}(I)+1-E_{00}(I)$.  Clearly,
$$ \psi_{n,s}|_{U_I(K)}: U_I(K)\simeq U_I(K),\;\;
\mu_I(\l )\mapsto \mu_I(\l ), $$ and $\psi_{n,s}(U_I(K)\ltimes
E_\infty (\mS_{CI}))=U_I(K)\ltimes E_\infty (L_{CI})$ for all
subsets $I$ with $|I|=s$. The subgroup of $(1+\ga_{n,s})^*$,
\begin{equation}\label{apns4}
\G_{n,s}:=\psi^{-1}_{n,s}(\bG_{n,s})={}^{set}\prod_{|I|=s}(U(K)\ltimes
E_\infty (\mS_{CI}))\cdot (1+\ga_{n,s+1})^*,
\end{equation}
is a normal subgroup as the  pre-image of a normal subgroup. We
added the upper script `set' to indicate that this is a product of
subgroups but not the direct product, in general. It is obvious
that $\psi_{n,s}(\G_{n,s})=\bG_{n,s}$ and $\G_{n,s}\subseteq
\Upsilon_{n,s}$. We will see that, in fact, $\G_{n,s}=
\Upsilon_{n,s}$ (Theorem \ref{30May9}).
  Let
$\D_{n,s}:=(1+\ga_{n,s})^*/\G_{n,s}$. The group homomorphism
$\psi_{n,s}$ (see (\ref{apns1}))  induces the group monomorphism
$$ \bpsi_{n,s}:\D_{n,s}\ra
\prod_{|I|=s}(1+\bgp_I)^*)/\bG_{n,s}\simeq \mX_{n,s}\simeq
\Z^{{n\choose s}(n-s)}.$$ This means that the group $\D_{n,s}$ is
a free abelian group of rank $\leq {n\choose s}(n-s)$. In fact,
the rank is equal to ${n\choose s+1}s$, see (\ref{Dns}).

For each subset $I$ with $|I|=s$, consider a free abelian group
$\mX_{CI}':= \bigoplus_{j\in CI}\Z (j,I)\simeq \Z^{n-s}$ where $\{
(j,I)\, | \, j\in CI\}$ is its free basis. Let
$$ \mX_{n,s}':=
\bigoplus_{|I|=s}\mX_{CI}'=\bigoplus_{|I|=s}\bigoplus_{j\in CI}\Z
(j,I)\simeq \Z^{{n\choose s}(n-s)}.$$ For each subset $I$,
consider the isomorphism of abelian groups
$$ \mX_{CI}\ra \mX_{CI}', \;\; \mu_I(x_j):= x_jE_{00}(I)+1-E_{00}(I)\mapsto (j,I).
$$
These isomorphisms yield the group isomorphism 
\begin{equation}\label{XXns}
\mX_{n,s}\ra \mX_{n,s}', \;\; \mu_I(x_j)\mapsto (j,I).
\end{equation}
Each element $a$ of the group $\mX_{n,s}$ is a unique product
$a=\prod_{|I|=s}\prod_{j\in CI}\mu_I(x_j)^{n(j,I)}$ where
$n(j,I)\in \Z$. Each element $a'$ of the group $\mX_{n,s}'$ is a
unique sum $a'=\sum_{|I|=s}\sum_{j\in CI}n(j,I)\cdot (j,I)$ where
$n(j,I)\in \Z$. The map (\ref{XXns}) sends $a$ to $a'$. To make
computations more readable we set $e_I:=E_{00}(I)$. Then
$e_Ie_J=e_{I\cup J}$.

$\noindent $

{\bf The current groups $\Theta_{n,s}$, $s=1,\ldots , n-1$}. The
current groups $\Theta_{n,s}$ are the most important subgroups of
the group $(1+\ga_n)^*$. They are finitely generated groups and
the generators are given explicitly. The adjective `current' comes
from the action of the generators on the monomial basis for  the
polynomial algebra $P_n$. When we visualize the algebra $P_n$ as a
liquid and the monomials $\{ x^\alpha \}$ are its atoms then the
action of the generators of the group $\Theta_{n,s}$ on the
monomials resembles  a current, see (\ref{curtij}). The generators
shift the liquid only on the faces of the positive cone
$\N^n\approx P_n$. The generators of the groups $\Theta_{n,s}$ are
units of the algebra $\mS_n$ but they are defined as a product of
two {\em non-units}. As a result the groups $\Theta_{n,s}$ capture
the most delicate phenomena about the structure and the properties
of the groups $\mS_n^*$ and $G_n$.

For each non-empty subset $I$ of $\{ 1,\ldots , n\}$ with
$s:=|I|<n$ and an element $i\in CI$, let
$$X(i,I):=\mu_I(x_i)= x_iE_{00}(I)+1-E_{00}(I)\;\; {\rm and}\;\;
Y(i,I):=\mu_I(y_i)= y_iE_{00}(I)+1-E_{00}(I).$$ Then
$Y(i,I)X(i,I)=1$, $\ker \, Y(i,I)=P_{C(I\cup i)}$, and $P_n=\im \,
X(i,I)\bigoplus P_{C(I\cup i)}$ where $P_{C(I\cup i)}:=
K[x_j]_{j\in C(I\cup i)}$. As an  element of the algebra
$\End_K(P_n)$, the map $X(i,I)$ is injective (but not bijective),
and the map $Y(i,I)$ is surjective (but not bijective).

$\noindent $

{\it Definition}. For each subset $J$ of $\{ 1,\ldots , n\}$ with
$|J|=s+1\geq 2$ and for two distinct elements $i$ and $j$ of the
set $J$,
$$ \th_{ij}(J):=Y(i,J\backslash i) X(j,J\backslash j)\in
(1+\gp_{J\backslash i}+\gp_{J\backslash j})^* \subseteq
(1+\ga_{n,s})^*.$$ The {\bf current group} $\Theta_{n,s}$ is the
subgroup of $(1+\ga_{n,s})^*$ generated by all the elements
$\th_{ij}(J)$ (for all the possible choices of $J$, $i,$ and $j$).

$\noindent $

In more detail, the element $\th_{ij}(J)$ belongs to the set
$1+\gp_{J\backslash i}+\gp_{J\backslash j}$ and
$\th_{ij}(J)^{-1}=\th_{ji}(J) \in 1+\gp_{J\backslash
i}+\gp_{J\backslash j}$. This follows from the action of the
element $\th_{ij}(J)$ on the monomial basis of the polynomial
algebra $P_n$, 
\begin{equation}\label{curtij}
\th_{ij}(J)*x^\alpha=\begin{cases}
x^\alpha & \text{if } \exists k \in J\backslash \{ i,j\}:\alpha_k\neq 0,\\
x^\alpha & \text{if } \forall k\in J\backslash \{ i,j\}:\alpha_k= 0, \alpha_i>0, \alpha_j>0,\\
x^{\alpha -e_i}& \text{if } \forall k\in J\backslash \{ i,j\}:\alpha_k= 0, \alpha_i>0, \alpha_j=0,\\
x^{\alpha +e_j}& \text{if } \forall k\in J\backslash \{ i,j\}:\alpha_k= 0, \alpha_i=0, \alpha_j\geq 0.\\
\end{cases}
\end{equation}
Alternatively, note that $\mu_{J\backslash j}(x_jy_j)=
\mu_{J\backslash j}(1-e_{\{ j\} }) = 1-e_{\{ j\} }e_{J\backslash
j} =1-e_J$ and (using (\ref{xyEij})) $\mu_{J\backslash
i}(y_i)e_J=(1+(y_i-1)e_{J\backslash i})e_J=e_J+(y_i-1)e_J=
e_J-e_J=0$, then
\begin{eqnarray*}
\th_{ij}(J)\th_{ji}(J)&= & \mu_{J\backslash i} (y_i)
\mu_{J\backslash j} (x_j)\cdot \mu_{J\backslash j}
(y_j)\mu_{J\backslash i } (x_i)=
\mu_{J\backslash i} (y_i)\cdot  \mu_{J\backslash j} (x_jy_j)\cdot \mu_{J\backslash i } (x_i)\\
 &=& \mu_{J\backslash i} (y_i)\cdot  (1-e_J)\cdot  \mu_{J\backslash i }
 (x_i)=\mu_{J\backslash i} (y_ix_i)= \mu_{J\backslash i } (1)=1.
\end{eqnarray*}
By symmetry, $\th_{ji}(J)\th_{ij}(J)=1$, i.e. 
\begin{equation}\label{tiji}
\th_{ij}(J)= \th_{ji}(J)^{-1}.
\end{equation}
 Therefore, the unit $\th_{ij}(I)$
is the product of an injective and surjective maps which are not
bijections.

Suppose that $i$, $j$, and $k$ are distinct elements of the set
$J$ (hence $|J|\geq 3$). Then 
\begin{equation}\label{tijjk}
\th_{ij}(J)\th_{jk}(J)=\th_{ik}(J).
\end{equation}
Indeed,
\begin{eqnarray*}
 \th_{ij}(J)\th_{jk}(J)&=& \mu_{J\backslash i}(y_i) \cdot \mu_{J\backslash j}(x_j)\mu_{J\backslash j}(y_j)\cdot \mu_{J\backslash k }(x_k) \\
 &=& \mu_{J\backslash i}(y_i) \cdot \mu_{J\backslash j}(x_jy_j)\cdot \mu_{J\backslash k }(x_k) \\
&=& \mu_{J\backslash i}(y_i) \cdot (1-e_J)\cdot \mu_{J\backslash k }(x_k) \\
&=& \mu_{J\backslash i}(y_i) \mu_{J\backslash k }(x_k)=
\th_{ik}(J).
\end{eqnarray*}
For each number $s=1, \ldots , n-1$, the free abelian group
$\mX_{n,s}'$ admits the decomposition
$\mX_{n,s}'=\bigoplus_{|J|=s+1}\bigoplus_{j\cup I=J}\Z (j,I)$, and
using it we define a character (a homomorphism) $\chi_J'$, for
each subset $J$ with $|J|=s+1$:
$$ \chi_J':\mX_{n,s}'\ra \Z, \;\; \sum_{|J'|=s+1}\sum_{j\cup I=J'}n_{j,I}(j,I)\mapsto \sum_{j\cup I=J}n_{j,I}.$$
Let $\max (J)$ be the maximal number  of the set $J$. The group
$\mX_{n,s}'$ is the direct sum 
\begin{equation}\label{XK1}
\mX_{n,s}'=\mK_{n,s}'\bigoplus \mY_{n,s}'
\end{equation}
of its free abelian subgroups,
\begin{eqnarray*}
 \mK_{n,s}'&=& \bigcap_{|J|=s+1}\ker (\chi_J')=\bigoplus_{|J|=s+1}\bigoplus_{j\in J\backslash
 \max (J)}\Z (-(\max (J), J\backslash \max (J))+(j,J\backslash j))\simeq \Z^{{n\choose s+1}s},  \\
 \mY_{n,s}'&=& \bigoplus_{|J|=s+1}\Z ((\max (J), J\backslash \max
 (J))\simeq \Z^{{n\choose s+1}}.
\end{eqnarray*}
Consider the group homomorphism $\psi_{n,s}':(1+\ga_{n,s})^*\ra
\mX_{n,s}'$ defined as the composition of the following group
homomorphisms:
$$\psi_{n,s}': (1+\ga_{n,s})^*\ra
(1+\ga_{n,s})^*/\G_{n,s}\stackrel{\psi_{n,s}}{\ra}
\prod_{|I|=s}(1+\bgp_I)^*/\bG_{n,s}\simeq \mX_{n,s}\simeq
\mX_{n,s}'.$$ Then 
\begin{equation}\label{ptijJ}
\psi_{n,s}'(\th_{ij}(J))=-(i,J\backslash i)+(j,J\backslash j).
\end{equation}
It follows that 
\begin{equation}\label{ptijJ1}
\psi_{n,s}'(\Theta_{n,s})=\mK_{n,s}',
\end{equation}
since, by (\ref{ptijJ}), $\psi_{n,s}'(\Theta_{n,s})\supseteq
\mK_{n,s}'$ (as the free basis for $\mK_{n,s}'$, introduced above,
belongs to the set $\psi_{n,s}'(\Theta_{n,s})$; again, by
(\ref{ptijJ}), $\psi_{n,s}'(\Theta_{n,s})\subseteq
\bigcap_{|J|=s+1}\ker (\chi_J')= \mK_{n,s}'$).

Let $H, H_1, \ldots , H_m$ be subsets (usually subgroups) of a
group $H'$.  We say that $H$ is the {\em  product} of $H_1, \ldots
, H_m$, and write $H={}^{set}\prod_{i=1}^mH_i=H_1\cdots H_m$, if
each element $h$ of $H$ is a  product $h=h_1\cdots h_m$ where
$h_i\in H_i$. We add the subscript `set' (sometime) in order to
distinguish it from the direct product of groups.
 We say that $H$ is the {\em exact product} of $H_1,
\ldots , H_m$, and write
$H={}^{exact}\prod_{i=1}^mH_i=H_1\times_{ex}\cdots
\times_{ex}H_m$, if each element $h$ of $H$ is a {\em unique}
product $h=h_1\cdots h_m$ where $h_i\in H_i$. The order in the
definition of the exact product is important.

The subgroup of $(1+\ga_{n,s})^*$ generated by the groups
$\Theta_{n,s}$ and $\G_{n,s}$ is equal to their product
$\Theta_{n,s}\G_{n,s}$, by the normality of $\G_{n,s}$. The
subgroup $\G_{n,s}$ of the group $\Theta_{n,s}\G_{n,s}$ is a
normal subgroup, hence the intersection $\Theta_{n,s}\cap
\G_{n,s}$ is a normal subgroup of $\Theta_{n,s}$.

\begin{lemma}\label{a24May9}
For each number $s=1, \ldots , n-1$, the group
$\Theta_{n,s}\G_{n,s}$ is the exact product
$$\Theta_{n,s}\G_{n,s}= {}^{exact}\prod_{|J|=s+1}\prod_{j\in J\backslash \max (J)}\langle
\th_{\max (J), j}(J)\rangle\cdot \G_{n,s},$$ i.e. each element
$a\in \Theta_{n,s}\G_{n,s}$ is a  unique product
$a=\prod_{|J|=s+1}\prod_{j\in J\backslash \max (J)} \th_{\max (J),
j}(J)^{n(j,J)}\cdot \g$ where $n(j,J)\in \Z$ and $\g \in
\G_{n,s}$. Moreover, the group $\Theta_{n,s}\G_{n,s}$ is the
semi-direct  product
$$\Theta_{n,s}\G_{n,s}= {}^{semi}\prod_{|J|=s+1}\prod_{j\in J\backslash \max (J)}\langle
\th_{\max (J), j}(J)\rangle\ltimes \G_{n,s},$$ where the order in
the double product is arbitrary.
\end{lemma}

{\it Proof}. The lemma follows at once  from (\ref{ptijJ1}) and
the fact that the elements $\psi_{n,s}'(\th_{\max (J),
j}(J))=-(\max (J), J\backslash \max (J))+(j,J\backslash j)$ form a
basis for the free abelian group $\mK_{n,s}'$.  $\Box $

$\noindent $

For each number $s=1, \ldots , n-1$, consider the subset of
$(1+\ga_{n,s})^*$, 
\begin{equation}\label{Tpns}
\Theta_{n,s}':= {}^{exact}\prod_{|J|=s+1}\prod_{j\in J\backslash
\max (J)}\langle \th_{\max (J), j}(J)\rangle,
\end{equation}
 which is the exact
product of cyclic groups (each of them is isomorphic to $\Z $)
since  each element $u$ of $\Theta_{n,s}'$ is a {\em unique}
product $u= \prod_{|J|=s+1}\prod_{j\in J\backslash \max
(J)}\th_{\max (J), j}(J)^{n(j,J)}$ where $n(j,J)\in \Z$ (Lemma
\ref{a24May9}).

By Lemma \ref{a24May9}, $\Theta_{n,s}/ \Theta_{n,s}\cap
\G_{n,s}\simeq \Theta_{n,s}\G_{n,s}/ \G_{n,s}\simeq
\mK_{n,s}'\simeq \Z^{{n\choose s+1}s}$, and so 
\begin{equation}\label{TnsG}
[\Theta_{n,s}, \Theta_{n,s}]\subseteq \G_{n,s}.
\end{equation}

The next theorem is the pennacle of finding  the explicit
generators for the groups $\mS_n^*$ and $G_n$.

\begin{theorem}\label{24May9}
$\psi_{n,s}'((1+\ga_{n,s})^*)=\psi_{n,s}'(\Theta_{n,s})$ for
$s=1,\ldots , n-1$.
\end{theorem}

{\it A rough sketch of the proof}. The proof is rather long, it is
given in Section \ref{SGNGEN}.  We use an  induction on $n$ (the
case $n=2$ was considered in  \cite{K1aut}), and then, for a fixed
$n$, we use a  second downward induction on $s=1, \ldots , n-1$
starting with $s=n-1$. The initial step $(n,s)=(n,n-1)$ is the
most difficult one. We spend entire Section \ref{GTNN1} to give
its proof. The remaining cases, using double induction, can be
deduced from the initial one (for different $n'$, i.e. when $n'$
runs from $1$ till $n$). The key idea in the proof of the case
$(n,n-1)$ is to use the Fredholm operators and their indices. Then
using well-known results on indices, some (new) results on the
Fredholm operators and their indices from \cite{K1aut} and their
generalizations obtained in Section \ref{GTNN1} we construct
several index maps (using various indices of the Fredholm
operators). The most difficult part is to prove that these maps
are well-defined (as their constructions are based on highly
non-unique decompositions). Then the proof follows from the
properties of these index maps. $\Box $

$\noindent $


\section{The groups $(1+\ga_{n, n-1})^*$ and $\Theta_{n,n-1}$} \label{GTNN1}

In this section, the group $(1+\ga_{n,n-1})^*$ is found (Corollary
\ref{a7May9}). We mentioned already in the Introduction that the
key idea in finding the group $G_n$ is to use indices of
operators. That is why we start this section with collecting known
results on indices and prove new ones. These results are used in
many  proofs that follow.

$\noindent $

{\bf The index $\ind$ of linear maps and its properties}. Let $\CC
$ be the family of all $K$-linear maps with finite dimensional
kernel and cokernel (such maps are called the {\em Fredholm linear
maps/operatorts}).  So, $\CC$ is the family of {\em Fredholm}
linear maps/operators. For vector spaces $V$ and $U$, let $\CC
(V,U)$ be the set of all the linear maps from $V$ to $U$ with
finite dimensional kernel and cokernel. So, $\CC =\bigcup_{V,U}\CC
(V,U)$ is the disjoint union.

$\noindent $

{\it Definition}. For a linear map $\v \in \CC$, the integer $
\ind (\v ) := \dim \, \ker (\v ) - \dim \, \coker (\v )$ is called
the {\em index} of the map $\v$.

$\noindent $

For vector spaces $V$ and $U$, let $\CC (V,U)_i:= \{ \v \in \CC
(V,U)\, | \, \ind (\v ) = i\}$. Then $\CC (V,U)=\bigcup_{i\in
\Z}\CC (V,U)_i$ is the disjoint union, and the set $\CC$ is the
disjoint union $\bigcup_{i\in \Z}\CC_i$ where $\CC_i :=\{ \v \in
\CC \, | \, \ind (\v ) = i\}$. When $V=U$, we write $\CC (V):= \CC
(V,V)$ and $\CC (V)_i:= \CC (V,V)_i$.

{\it Example}. Note that $\mS_1\subset \End_K(P_1)$. The map
$x^i\in \End_K(P_1)$  is an injection with
$P_1=(\bigoplus_{j=0}^{i-1}Kx^j)\bigoplus \im (x^i)$; and the map
$y^i\in \End_K(P_1)$  is a surjection with $\ker (y^i)
=\bigoplus_{j=0}^{i-1}Kx^j$. Hence  
\begin{equation}\label{indxy}
\ind (x^i)= -i\;\; {\rm and }\;\; \ind (y^i)= i, \;\; i\geq 1.
\end{equation}

Lemma \ref{b8Feb9} shows that $\CC$ is a multiplicative semigroup
with zero element (if the composition of two elements of $C$ is
undefined we set their product to be zero). The next two lemmas
are well known.

\begin{lemma}\label{b8Feb9}
Let $\psi : M\ra N $ and $\v : N\ra L$ be $K$-linear maps. If two
of the following three maps: $\psi$, $\v$,  and $\v \psi$, belong
to the set $\CC$ then so does the third; and in this case, $$ \ind
(\v\psi ) = \ind (\v ) + \ind (\psi ).$$
\end{lemma}
 By lemma \ref{b8Feb9}, $\CC (N,L)_i\CC (M,N)_j\subseteq \CC
 (M,L)_{i+j}$ for all $i,j\in \Z$.
\begin{lemma}\label{a8Feb9}
Let
$$
\xymatrix{0\ar[r] & V_1\ar[r]\ar[d]^{\v_1}  & V_2 \ar[r]\ar[d]^{\v_2} & V_3 \ar[r]\ar[d]^{\v_3 } & 0 \\
0\ar[r] & U_1\ar[r]  & U_2\ar[r] & U_3 \ar[r] & 0 }
$$
be a commutative diagram of $K$-linear maps with exact rows.
Suppose that $\v_1, \v_2, \v_2\in \CC$. Then
$$ \ind (\v_2) = \ind (\v_1)+\ind (\v_3).$$
\end{lemma}

Let $V$ and $U$ be vector spaces. Define $\CI (V,U) :=\{ \v \in
\Hom_K(V,U)\, | \, \dim \, \im (\v ) <\infty \}$, and when $V=U$
we write $\CI (V):= \CI (V,V)$.

\begin{theorem}\label{b23Apr9}
\cite{K1aut} Let $V$ and $U$ be vector spaces. Then $\CC
(V,U)_i+\CI (V, U) = \CC (V,U)_i$ for all $i\in \Z$.
\end{theorem}

\begin{lemma}\label{a21Apr9}
\cite{K1aut} Let $V$ and $V'$ be vector spaces, and  $\v :V\ra V'$
be a linear map such that the vector spaces $\ker (\v )$ and
$\coker (\v )$
 are isomorphic. Fix subspaces $U\subseteq V$ and $W\subseteq V'$
 such that $V=\ker (\v ) \bigoplus U$ and $V'= W\bigoplus \im (\v
 )$  and fix an isomorphism $f: \ker (\v ) \ra W$ (this is
 possible since $\ker (\v ) \simeq \coker (\v ) \simeq W$) and
 extend it to a linear map $f: V\ra V'$ by setting $f(U)=0$. Then
 the map $\v +f:V\ra V'$ is an isomorphism.
\end{lemma}

\begin{corollary}\label{a23Apr9}
\begin{enumerate}
\item $1+F_n\subseteq \CC (P_n)_0$. \item $\mS_n^* +F_n\subseteq
\CC (P_n)_0$.
\end{enumerate}
\end{corollary}

{\it Proof}. Both statements follows from Theorem \ref{b23Apr9}:
$\mS_n^*\in \CC (P_n)_0$ and $F_n\in \CI (P_n)$, but we give short
independent proofs (that do not use Theorem \ref{b23Apr9}).

 1. Since $1+F_n\simeq 1+M_\infty (K)$, statement 1 is obvious.

2. Let $u\in \mS_n^*$ and $f\in F_n$. Then $u^{-1}f\in F_n$. By
statement 1, the element $1+u^{-1}f\in \CC (P_n)_0$. Since $u\in
\CC (P_n)_0$, we have $u+f= u(1+u^{-1}f)\in \CC (P_n)_0$, by Lemma
\ref{b8Feb9}.  $\Box $

$\noindent $

{\bf The subgroup $\Theta_{n,n-1}$ of $(1+\ga_{n,n-1})^*$ for
$n\geq 2$}. For each pair of indices $i\neq j$,  the element
$$\th_{ij} :=\th_{ij}(\{ 1,\ldots , n\}) := (1+(y_i-1)\prod_{k\neq
i}E_{00}(k))\cdot (1+(x_j-1)\prod_{l\neq j}E_{00}(l))\in
(1+\ga_{n,n-1})^*$$ is a unit and 
\begin{equation}\label{thijm1}
\th^{-1}_{ij}= (1+(y_j-1)\prod_{l\neq j}E_{00}(l))\cdot
(1+(x_i-1)\prod_{k\neq i}E_{00}(k))\in (1+\ga_{n,n-1})^*,
\end{equation}
i.e. $\th_{ij}^{-1}=\th_{ji}$.  This is obvious since
$$\th_{ij} *x^\alpha = \begin{cases}
x_i^{\alpha_i-1}& \text{if } \alpha_i>0, \forall k\neq i: \alpha_k=0, \\
x_j^{\alpha_j+1}& \text{if }  \alpha_j\geq 0, \forall l\neq j: \alpha_l=0,  \\
x^\alpha& \text{otherwise},\\
\end{cases}
\;\;\; {\rm and}\;\;\; \th^{-1}_{ij} *x^\alpha = \begin{cases}
x_i^{\alpha_i+1}& \text{if }  \alpha_i\geq 0, \forall k\neq i: \alpha_k=0, \\
x_j^{\alpha_j-1}& \text{if } \alpha_j> 0,  \forall l\neq j: \alpha_l=0. \\
x^\alpha& \text{otherwise}.\\
\end{cases}
$$
Using the above action of the elements $\th_{ij}$ on the monomial
basis for the polynomial algebra $P_n$, it is easy to show that
the elements $\th_{ij}$ {\em commute} modulo $(1+F_n)^*$; ,
$\th_{jk}\th_{ij} \equiv \th_{i,k} \mod (1+F_n)^*$ for all
distinct elements  $i$, $j$, and $k$; and $\th_{ij}*1=x_j^m$ for
all $m\geq 1$. Recall that $\Theta_{n,n-1}$ is the subgroup of
$(1+\ga_{n,n-1})^*$ generated by the elements $\th_{ij}$. It
follows from
$$ (1+(y_i-1)\prod_{k\neq i}E_{00}(k))*x^\alpha = \begin{cases}
x_i^{\alpha_i-1}& \text{if } \alpha_i>0, \forall k\neq i: \alpha_k=0, \\
0& \text{if } \alpha=0,  \\
x^\alpha& \text{otherwise},\\
\end{cases}$$
that the map $1+(y_i-1)\prod_{k\neq i}E_{00}(k)\in \End_K(P_n)$ is
a surjection with kernel equal to $K$, and so 
\begin{equation}\label{nindEy}
\ind (1+(y_i-1)\prod_{k\neq i}E_{00}(k))=1.
\end{equation}
 Similarly, it
follows from
$$ (1+(x_j-1)\prod_{l\neq
j}E_{00}(l))*x^\alpha = \begin{cases}
x_j^{\alpha_j+1}& \text{if } \forall l\neq j: \alpha_l=0,  \\
x^\alpha& \text{otherwise},\\
\end{cases}$$
that the map $1+(x_j-1)\prod_{l\neq j}E_{00}(l)\in \End_K(P_n)$ is
an injection such that  $P_n= K\bigoplus \im
(1+(x_j-1)\prod_{l\neq j}E_{00}(l))$, and so 
\begin{equation}\label{nindEy1}
\ind (1+(x_j-1)\prod_{l\neq j}E_{00}(l))=-1.
\end{equation}
 We see that the unit $\th_{ij}$ of the algebra
$\mS_n$ is the product of two non-units having nonzero indices of
opposite sign  (note that $\ind (\th_{ij} )=0$, and so  the sum of
the two indices is equal to zero). Lemma \ref{a24Apr9} shows that
this is a general phenomenon, and so the group $(1+\ga_{n,n-1})^*$
is a `complicated' group in the sense that in producing units
non-units are involved.

\begin{lemma}\label{a24Apr9}
Let $u=1+\sum_{i=1}^na_i\in (1+\ga_{n, n-1})^*$ where $a_i\in
\gp_{Ci}$. Then
\begin{enumerate}
\item  $1+a_i\in \CC (P_n)$ for all $i=1, \ldots, n$; and
$\sum_{i=1}^n\ind (1+a_i) =0$. \item If $u=1+\sum_{i=1}^na_i'$
where $a_i'\in \gp_{Ci}$ then $\ind (1+a_i) = \ind (1+a_i')$ for
all $i=1, \ldots, n$.
\end{enumerate}
\end{lemma}

{\it Proof}. 1. Since $a_i\in \gp_{Ci}$ for all $i$, we have
$a_ia_j\in F_n$ provided $i\neq j$.
 It follows that the elements $f:= u-(1+a_1)(1+a_2)\cdots (1+a_n)$
 and $f':= u-(1+a_2)\cdots (1+a_n)(1+a_1)$ belong to the ideal
 $F_n$. By Corollary \ref{a23Apr9}.(2), $u-f, u-f'\in \CC
 (P_n)_0$. Then, it follows from the equalities  $u-f= (1+a_1)(1+a_2)\cdots (1+a_n)$
 and $ u-f'=(1+a_2)\cdots (1+a_n)(1+a_1)$ that
$$
 \im (1+a_1) \supseteq  \im ( u-f) \; \; {\rm and}\; \; \ker (1+a_1) \subseteq  \ker ( u-f').
$$
This means that $1+a_1\in \CC (P_n)$. By symmetry, $1+a_i\in \CC
(P_n)$ for all $i$.
 By Corollary
\ref{a23Apr9}.(2) and Lemma \ref{b8Feb9}, $$0=\ind (u)= \ind (u-f)
= \ind (1+a_1)\cdots (1+a_n) = \sum_{i=1}^n\ind (1+a_i).$$ 2. For
each number  $i$, $f_i:= a_i'-a_i= -\sum_{j\neq i}(a_j'-a_j)\in
\gp_{Ci}\bigcap  \gp_i= \bigcap_{j=1}^n \gp_j= F_n$. Since
$F_n\subseteq \CI (P_n)$, we see that $\ind (1+a_i') = \ind
(1+a_i+f_i) = \ind (1+a_i)$, by Theorem \ref{b23Apr9}. $\Box $

$\noindent $

By Lemma \ref{a24Apr9}, for each number  $i=1, \ldots , n$, there
is a well-defined map,
\begin{equation}\label{indi}
\ind_i: (1+\ga_{n,n-1})^*\ra \Z , \;\; u=1+\sum_{i=1}^n a_i\mapsto
\ind (1+a_i),
\end{equation}
(where $a_i\in \gp_{Ci}$) which is a group homomorphism:
\begin{eqnarray*}
 \ind_i(uu') &= & \ind_i((1+\sum_{i=1}^n a_i)(1+\sum_{j=1}^n a_j') )= \ind(1+a_i+a_i'+a_ia_i')     \\
 &=& \ind ((1+a_i)(1+a_i') )= \ind (1+a_i)+\ind(1+a_i')\\
 & =& \ind_i(u) +\ind_i(u'),
\end{eqnarray*}
since $a_ja_j'\in \gp_{Cj}$ for all $j$.  Let $\CK_{n,n-1}$ be the
kernel of the group epimorphism
$$ \bigoplus_{i=1}^{n-1}\ind_i : (1+\ga_{n,n-1})^*\ra
\Z^{n-1}=\bigoplus_{i=1}^{n-1} \Z e_i, \;\; 1+\sum_{i=1}^n
a_i\mapsto \sum_{i=1}^{n-1}\ind (1+a_i)\cdot e_i, $$ where $a_i\in
\gp_{Ci}$ for $i=1, \ldots , n$. The restriction of the
epimorphism to the subset
$\Theta_{n,n-1}':={}^{exact}\prod_{j=1}^{n-1}\langle
\th_{n,j}\rangle $ is a bijection since (by (\ref{nindEy}) and
(\ref{nindEy1}))
$$\bigoplus_{i=1}^{n-1}\ind_i (\th_{j,j+1})= \begin{cases}
e_j-e_{j+1}& \text{if } j<n-1,\\
e_{n-1}& \text{if } j=n-1.\\
\end{cases}$$
Therefore, 
\begin{equation}\label{1aat}
(1+\ga_{n,n-1})^* = {}^{exact} \Theta_{n,n-1}' \cdot \CK_{n,n-1},
\;\; \CK_{n,n-1}=\bigcap_{i=1}^n \ker ( \ind_i),
\end{equation}
by Lemma \ref{a24Apr9}.(1). So, $\CK_{n,n-1}$ is the normal
subgroup of the group $(1+\ga_{n,n-1})^*$, $\Theta_{n,n-1}'\cap
\CK_{n,n-1}=\{ 1\}$, and each element $u$ of the group
$(1+\ga_{n,s})^*$ is a unique product $vw$ for some elements $v\in
\Theta_{n,n-1}'$ and $w\in \CK_{n,n-1}$.  The subgroups
$(1+\gp_{Ci})^*$, $i=1, \ldots , n$ of the group
$(1+\ga_{n,n-1})^*$ or $(1+\ga_n)^*$ are normal, and
$(1+\gp_{Ci})^*\bigcap (1+\gp_{Cj})^*=(1+F_n)^*$ for all $i\neq
j$. The product $\prod_{i=1}^n (1+\gp_{Ci})^*:= \{ u_1\cdots u_n\,
| \, u_i\in (1+\gp_{Ci})^*,   i=1, \ldots ,n\}$ is a normal
subgroup of $(1+\ga_{n,n-1})^*$ and $(1+\ga_n)^*$. In fact, the
order in the product can be arbitrary (by normality). Clearly,
$\prod_{i=1}^n (1+\gp_{Ci})^*\subseteq \CK_{n,n-1}$. In fact, the
equality holds as the next proposition shows.

\begin{proposition}\label{b24Apr9}
\begin{enumerate}
\item $\CK_{n,n-1}=\prod_{i=1}^n (1+\gp_{Ci})^*$.\item
$(1+\ga_{n,n-1})^*= {}^{exact} \Theta_{n,n-1}'\cdot (\prod_{i=1}^n
(1+\gp_{Ci})^*)=\langle \th_{n,1}\rangle\ltimes \cdots \ltimes
\langle \th_{n,n-1}\rangle \ltimes (\prod_{i=1}^n
(1+\gp_{Ci})^*)$.
\end{enumerate}
\end{proposition}

{\it Proof}. 1. It suffices to show that each element
$u=1+\sum_{i=1}^n a_i$ (where $a_i\in \gp_{Ci}$) of the group
$\CK_{n,n-1}$ is a product $u_1\cdots u_n$ of  some elements
$u_i\in (1+\gp_{Ci})^*$. By Lemma \ref{a24Apr9}, $1+a_1\in \CC
(P_n)_0$ since $u\in \CK_{n,n-1}$. Fix a subspace, say $W$, of
$P_n$ such that $P_n= \ker (1+a_1) \bigoplus W$ and
$W=\bigoplus_{\alpha \in I} Kx^\alpha$ where $I$ is a subset of
$\N^n$. By Lemma \ref{a21Apr9}, we can find an element $f_1\in
F_n$ (since $\dim \, \ker (1+a_1) <\infty $, $W$ has a monomial
basis, and $f_1(W)=0$) such that $u_1:= 1+a_1+f_1\in \Aut_K(P_n)$.
We claim that $u_1\in (1+\gp_{C1})^*$. It is a subtle point since
{\em not all} elements of the algebra $\mS_n$ that are invertible
linear maps in $P_n$ are invertible in $\mS_n$, i.e. $\mS_n^*
\subsetneqq \mS_n\cap \Aut_K(P_n)$ but $$(1+F_n)^* = (1+F_n)\cap
\Aut_K(P_n),\;\; \cite{shrekaut}. $$ The main idea in the proof of
the claim is to use this equality. Similarly, for each $i\geq 2$,
we can find an element $f_i\in F_n$ such that $v_i:= 1+a_i+f_i\in
\Aut_K(P_n)$. Then $v:=v_2\cdots v_n\in \Aut_K(P_n)$,  $u=
u_1v+g_1$ and $u= vu_1+g_2$ for some elements $g_i\in F_n$. Hence,
$$ u_1vu^{-1}=1-g_1u^{-1}\;\; {\rm and }\;\;\ u^{-1} vu_1= 1-u^{-1}
g_2, $$ and so $1-g_1u^{-1}, 1-u^{-1} g_2\in (1+F_n) \cap
\Aut_K(P_n) = (1+F_n)^*$. It follows that $u_1^{-1} = vu^{-1}
(1-g_1u^{-1})^{-1}\in (1+\gp_{C1})^*$ since
$$1\equiv  1-g_1u^{-1}
\equiv u_1vu^{-1} \equiv vu^{-1}\mod \gp_{C1}.$$ This proves the
claim. Clearly, $u_2':= v+u_1^{-1}g_1\in 1+\sum_{j=2}^n\gp_{Cj}$.
Then, it follows from the equality $ u = u_1v+g_1 =
u_1(v+u_1^{-1}g_1)= u_1u_2'$ that $u_2'= u_1^{-1} u\in
(1+\sum_{j=2}^n\gp_{Cj})^*$. Repeating the same argument for the
element $u_2'$ we find an element $ u_2\in (1+\gp_{C2})^*$ such
that $u_3':= u_2^{-1}u_2'\in (1+\sum_{j=3}^n \gp_{Cj})^*$.
Repeating the same argument again and again (or use induction) we
find elements $u_i\in (1+\gp_{Ci})^*$ and elements $u_i'\in
(1+\sum_{j=i+1}^n\gp_{Cj})^*$ such that
$u_i'=u_{i-1}^{-1}u_{i-1}'$, hence $u=u_1u_2'=u_1u_2u_3'=\cdots
=u_1u_2\cdots u_n$, as required.

2. Statement 2 follows from statement 1 and (\ref{1aat}).   $\Box
$

$\noindent $

 For each number $i=1,\ldots , n$, the group of units
of the monoid $1+\gp_{Ci } = 1+\mS_1(i)\bigotimes
\bigotimes_{j\neq i}F(j)\simeq 1+M_\infty (\mS_1(i))$ is equal to
$(1+\gp_{Ci})^* \simeq  \GL_\infty (\mS_1(i))$. It contains the
semi-direct product $U_{Ci}(K)\ltimes E_\infty (\mS_1(i))$ of its
two subgroups, where
$$U_{Ci}(K):=\{ \l \prod_{j\neq i}E_{00}(j) +1-\prod_{j\neq i}E_{00}(j)\, | \, \l \in K^*\} \simeq
K^*,$$
 and the group $E_\infty (\mS_1(i))$ is generated by all the
 elementary matrices $1+aE_{kl}(Ci)$ where $k,l\in \N^{n-1}$,
 $k\neq l$, $E_{kl}(Ci):=\prod_{j\neq i}E_{k_jl_j}(j)$,   and $a\in
 \mS_1(i)$. We will see (Proposition \ref{c24Apr9}) that the
 group  $(1+\gp_{Ci})^*$ coincides with the semi-direct product.

The set $F_n$ is an ideal of the algebra $K+\gp_{Ci}=
K(1+\gp_{Ci})$ which is a subalgebra  of the algebra $\mS_n$, and
$(K+\gp_{Ci}) / F_n= K(1+\gp_{Ci}/ F_n) \simeq K(1+M_\infty
(L_{i}))$ where $L_{i}:= K[x_{i}, x_{i}^{-1}]\simeq \mS_1(i)/F(i)$
is the Laurent polynomial algebra. The algebra $L_{i}$ is a
 Euclidean  domain, hence $\GL_\infty (L_{i}) =
U(L_{i})\ltimes E_\infty (L_{i})$ where
$$U(L_{i}):=\{ a \prod_{j\neq i}
E_{00}(j) +1-\prod_{j\neq i}E_{00}(j)\, | \, a \in L_{i}^*\}
\simeq L_{i}^*=K^*\times \{ x_{i}^m \, | \, m\in \Z\}$$ and
$E_\infty (L_{i})$ is the subgroup of $\GL_\infty (L_{i})$
generated by all the elementary matrices.

The group of units of the algebra $(K+\gp_{Ci}) / F_n$ is equal to
$K^*\times \GL_\infty (L_{i})= K^* \times ( U(L_{i})\ltimes
E_\infty (L_{i}))$. The algebra epimorphism $\psi_{Ci} :
K+\gp_{Ci}\ra (K+\gp_{Ci})/F_n$, $a\mapsto a+F_n$, induces the
exact sequence of groups, 
\begin{equation}\label{1piex}
1\ra (1+F_n)^*\ra (1+\gp_{Ci})^*\stackrel{\psi_{Ci} }{\ra}
\GL_\infty (L_{i})=U(L_{i})\ltimes E_\infty (L_{i}),
\end{equation}
which yields the short exact sequence of groups, 
\begin{equation}\label{2piex}
1\ra (1+F_n)^*\ra U_{Ci}(K)\ltimes E_\infty (\mS_1(i))\ra
U(K)\ltimes E_\infty (L_{i})\ra 1.
\end{equation}
Recall that $U_{Ci}(K)\ltimes E_\infty (\mS_1(i))\subseteq
(1+\gp_{Ci})^*$.  In fact, the equality holds.

\begin{proposition}\label{c24Apr9}
$(1+\gp_{Ci})^* =U_{Ci}(K)\ltimes E_\infty (\mS_1(i))$ for all
$i=1,\ldots , n$.
\end{proposition}

{\it Proof}. In view of the exact sequences (\ref{1piex}) and
(\ref{2piex}), it suffices to show that the image of the map
$\psi_{Ci}$ in (\ref{1piex}) is equal to $U(K)\ltimes E_\infty
(L_{i}))$. Since
$$U(L_{i}) = U(K)\times \{x_{i}^m \prod_{j\neq i}E_{00}(j)  +1- \prod_{j\neq i}E_{00}(j)  \, | \,
m\in \Z\},$$
 this is equivalent  to show that that if $\psi_{Ci} (u) =x_{i}^m \prod_{j\neq i}
E_{00}(j) +1-\prod_{j\neq i}E_{00}(j)$ for some element $u\in
(1+\gp_{Ci})^*$ and an integer $m\in \Z$ then $m=0$. Let $u(m) :=
v_{i}(m)\prod_{j\neq i}E_{00}(j)+1-\prod_{j\neq i}E_{00}(j)$ where
$v_{i}(m):=
\begin{cases}
x_{i}^m& \text{if } m\geq 0,\\
y_{i}^{|m|} & \text{if } m<0.\\
\end{cases} $ Then $u(m) \in 1+\gp_{Ci}$ and $\psi_{Ci} (u(m)) =
\psi_{Ci} (u)$. Hence, $u(m) = u+f_m$ for some element $f_m\in
F_n$. Note that $$u(m) =\begin{cases}
u(1)^m& \text{if } m\geq 0,\\
u(-1)^{|m|} & \text{if } m<0, \\
\end{cases} $$ and, by (\ref{nindEy}) and (\ref{nindEy1}), $\ind
(u(m))=-m$. By Corollary \ref{a23Apr9}.(2).
$$ 0=\ind (u) = \ind (u+f_m) = \ind (u(m)) = -m,$$
and so $m=0$, as required. $\Box $

$\noindent $

Combining Proposition \ref{b24Apr9}.(1) and Proposition
\ref{c24Apr9}, we have the next corollary.

\begin{corollary}\label{a7May9}
$(1+\ga_{n,n-1})^* = \Theta_{n,n-1}'\times_{ex}
({}^{set}\prod_{i=1}^n (1+\gp_{Ci})^*)\simeq
\Theta_{n,n-1}'\times_{ex} ({}^{set}\prod_{i=1}^n U_{Ci}(K)\ltimes
E_\infty (\mS_1(i)))\simeq \langle \th_{n,1}\rangle\ltimes \cdots
\ltimes \langle \th_{n,n-1}\rangle \ltimes ({}^{set}\prod_{i=1}^n
U_{Ci}(K)\ltimes E_\infty (\mS_1(i)))$.
\end{corollary}

Using Corollary \ref{a7May9}, we can write down explicit
generators for the group $(1+\ga_{n,n-1})^*$, see Theorem
\ref{25May9} where explicit generators are given for all the
groups $(1+\ga_{n,s})^*$.


\section{The structure of the  groups $\mS_n^*$, $G_n$ and their generators}\label{SGNGEN}

In this section, a proof of Theorem \ref{24May9} is given, the
groups $\mS_n^*$, $(1+\ga_n)^*$, and $G_n$ and their generators
are found explicitly (Theorem \ref{A10May9}, Theorem
\ref{B10May9}, Theorem \ref{25May9}, and Theorem \ref{A25May9}).

$\noindent $

{\bf Proof of Theorem \ref{24May9}.} To prove the theorem we use
induction on $n$. The inial step when $n=2$ follows from Corollary
\ref{a7May9} as in this case there is only one option,
$(n,s)=(2,1)$. So, let $n>2$, and suppose that the theorem holds
for all pairs $(n',s')$, $s'=1, \ldots , n'-1$ such that $n'<n$.
For the number $n$, we use a second downward  induction on $s=1,
\ldots , n-1$ starting with $s=n-1$. In this case, i.e.
$(n,s)=(n,n-1)$, the theorem holds as  it follows from Corollary
\ref{a7May9}. So, let $s<n-1$, and  suppose that the statement is
true for all pairs $(n, s')$ with $s'=s+1, \ldots , n-1$. For each
 number  $i=1, \ldots , n$, the algebra $\mS_{Ci}\t K(x_i)$ is isomorphic
to the algebra $\mS_{n-1}$ but over the field $K(x_i)$ of rational
functions. By the  induction on $n$, the theorem holds for the
algebra $\mS_{Ci}\t K(x_i)$. In order to stress that we consider
the algebra $\mS_{Ci}$ over the field $K(x_i)$ rather than $K$  we
add the subscript `$Ci$' to all the notations introduced for the
algebra $\mS_{Ci}$ but over the field $K$. For example, $\ga_{n-1,
s,Ci}\t K(x_i)$ stands for the ideal $\ga_{n-1,s}$ of the algebra
$\mS_{Ci}$ but over the field $K(x_i)$, etc.

For each number $i=1, \ldots , n$ and for each number $s=1, \ldots
, n-2$, the composition of the two algebra homomorphisms
$$ \mS_n\ra \mS_n/\gp_i\simeq \mS_{Ci}\t K[x_i, x_i^{-1}]\ra \mS_{Ci}\t
K(x_i)$$ induces the group homomorphism $(1+\ga_{n,s})^*\ra
(1+\ga_{n-1, s,Ci}\t K(x_i))^*$. This homomorphism yields the
commutative diagram where all the maps are obvious (and natural):

$$
\xymatrix{ (1+\ga_{n,s})^*\ar[r]\ar[d] & \frac{(1+\ga_{n,s})^*}{\G_{n,s}}\ar[r]^{\bpsi_{n,s}}\ar[d] & \frac{\prod_{|I|=s}(1+\bgp_I)^*}{\bG_{n,s}} \ar[d]\simeq \mX_{n,s}\simeq \mX_{n,s}'  \ar[d]^{\v_{n,s,i}} & \\
 (1+\ga_{n-1,s,Ci}\t K(x_i))^*\ar[r]& \frac{(1+\ga_{n-1,s,Ci}\t K(x_i))^*}{\G_{n-1,s,Ci}}\ar[r]^{\bpsi_{n-1,s, Ci}}& \; \;\;\;\;\;\;\;\;\; \frac{\prod'(1+\bgp_I\t K(x_i))^*}{\bG_{n-1,s, Ci}} \simeq \mX_{n-1,s, Ci}\simeq \mX_{n-1,s, Ci}'  &}
$$
where $\prod':= \prod_{ \{ I: |I|=s, i\not\in I \} }$ and the map
$\v_{n,s,i}: \mX_{n,s}'\ra \mX_{n-1,s,Ci}'$ is given by the rule
$$\v_{n,s,i}((j,I))=\begin{cases}
(j,I)& \text{if }i\notin I\cup j,\\
0& \text{otherwise}.\\
\end{cases} $$ This is obvious. By the induction on $n$, we have the
equality $\psi_{n-1, s, Ci}'((1+\ga_{n-1, s,Ci}\t
K(x_i))^*)=\psi_{n-1,s,Ci}'(\Theta_{n-1, s,Ci})$ for each $s=1,
\ldots , n-2$. Then, by the commutative diagram above,
\begin{equation}\label{pps1}
\v_{n,s,i}\psi_{n,s}'((1+\ga_{n,s})^*)\subseteq \psi_{n-1, s,
Ci}'((1+\ga_{n-1, s,Ci}\t K(x_i))^*)=\psi_{n-1,s,Ci}'(\Theta_{n-1,
s,Ci}).
\end{equation}
It follows from the definition of the map $\v_{n,s,i}$ that
\begin{equation}\label{pps2}
\v_{n,s,i}(\mY_{n,s}')\subseteq \mY_{n-1,s,Ci}'.
\end{equation}
Summarizing, for each $i=1, \ldots , n$, by (\ref{XK1}) and
(\ref{ptijJ1}), there is the map
$$ \v_{n,s,i}:\mX_{n,s}'=\psi_{n,s}'(\Theta_{n,s})\bigoplus
\mY_{n,s}'\ra
\mX_{n-1,s,Ci}'=\psi_{n-1,s,Ci}(\Theta_{n-1,s,Ci})\bigoplus\mY_{n-1,s,Ci}'$$
satisfying (\ref{pps1}) and (\ref{pps2}).  The group homomorphism
$\v_{n,s}:=\prod_{i=1}^n\v_{n,s,i}:\mX_{n,s}'\ra
\prod_{i=1}^n\mX_{n-1, s,Ci}'$ is a monomorphism since it has
trivial kernel: $\ker (\v_{n,s})=\bigoplus_{i\in I\cup j}\Z (j,I)$
where the pairs $(j,I)$ in the direct sum are such that $i\in
I\cup j$ for all $i=1, \ldots , n$ (see the definition of the map
$\v_{n,s,i}$), i.e. $I\cup j= \{ 1, \ldots , n\}$ but the number
of elements in the set $I\cup j$ is $s+1<n-1+1=n$, a
contradiction. This means that $\ker (\v_{n,s})=0$. Let $u\in
(1+\ga_{n,s})^*$. Then $\psi_{n,s}'(u)=a+b$ for unique elements
$a\in \psi_{n,s}'(\Theta_{n,s})$ and $b\in \mY_{n,s}'$. By
(\ref{pps1}) and (\ref{pps2}), $\v_{n,s,i}(b)=0$ for all $i=1,
\ldots , n$, i.e. $\v_{n,s}(b)=0$, and so $b=0$ since the map
$\v_{n,s}$ is a monomorphism. This proves that
$\psi_{n,s}'((1+\ga_{n,s})^*)=\psi_{n,s}(\Theta_{n,s})$. By
induction, the theorem holds. The proof of Theorem \ref{24May9} is
complete. $\Box$

$\noindent $

For each number $s=1, \ldots , n-1$, consider the following
subsets of the group $(1+\ga_{n,s})^*$, 
\begin{equation}\label{mEns}
\mE_{n,s}:=\prod_{|I|=s}U_I(K)\ltimes E_\infty (\mS_{CI}) \;\;
{\rm and}\;\;  \mP_{n,s}:=\prod_{|I|=s}(1+\gp_i)^*,
\end{equation}
 the products
of subgroups of $(1+\ga_{n,s})^*$ in arbitrary order which is
fixed for each $s$.

\begin{theorem}\label{A10May9}
\begin{enumerate}
\item $(1+\ga_n)^* = \Theta_{n,1}\G_{n,1} = \Theta_{n,1} \mE_{n,1}
\Theta_{n,2} \mE_{n,2}\cdots \Theta_{n,n-1} \mE_{n,n-1} $.
Moreover,  for $s=1, \ldots , n-1$,  $(1+\ga_{n,s})^* =
\Theta_{n,s}\G_{n,s} = \Theta_{n,s} \mE_{n,s} \Theta_{n,s+1}
\mE_{n,s+1}\cdots \Theta_{n,n-1} \mE_{n,n-1} $. \item $(1+\ga_n)^*
= \Theta_{n,1}\Upsilon_{n,1}=
 \Theta_{n,1} \mP_{n,1} \Theta_{n,2}
\mP_{n,2}\cdots \Theta_{n,n-1} \mP_{n,n-1} $. Moreover,  for $s=1,
\ldots , n-1$, $(1+\ga_{n,s})^* = \Theta_{n,s}\Upsilon_{n,s}=
 \Theta_{n,s} \mP_{n,s} \Theta_{n,s+1}
\mP_{n,s+1}\cdots \Theta_{n,n-1} \mP_{n,n-1} $.
\end{enumerate}
\end{theorem}

{\it Proof}. 1. By Theorem \ref{24May9} and Corollary
\ref{a7May9},
\begin{eqnarray*}
 (1+\ga_{n,s})^*&=&\Theta_{n,s}\G_{n,s}=\Theta_{n,s}\prod_{|I|=s}U_I(K)\ltimes
E_\infty (\mS_{CI})\cdot (1+\ga_{n,s-1})^*=\Theta_{n,s}
\mE_{n,s}(1+\ga_{n,s-1})^*\\
&=&\Theta_{n,s} \mE_{n,s}\Theta_{n,s-1}
\mE_{n,s-1}(1+\ga_{n,s-2})^*= \Theta_{n,s} \mE_{n,s}\cdots
\Theta_{n,n-2} \mE_{n,n-2}(1+\ga_{n,n-1})^*\\
&=& \Theta_{n,s} \mE_{n,s} \Theta_{n,s+1} \mE_{n,s+1}\cdots
\Theta_{n,n-1} \mE_{n,n-1}.
\end{eqnarray*}

2. Since $(1+\ga_{n,s})^*=\Theta_{n,s}\G_{n,s}\subseteq
\Theta_{n,s}\Upsilon_{n,s}\subseteq (1+\ga_{n,s})^*$, we see that
\begin{eqnarray*}
(1+\ga_{n,s})^*&=&\Theta_{n,s}\Upsilon_{n,s}=
\Theta_{n,s}\prod_{|I|=s}(1+\gp_I)^*\cdot (1+\ga_{n,s-1})^*=\Theta_{n,s} \mP_{n,s}(1+\ga_{n,s-1})^*\\
&=&\Theta_{n,s} \mP_{n,s}\Theta_{n,s-1}
\mP_{n,s-1}(1+\ga_{n,s-2})^*= \Theta_{n,s} \mP_{n,s}\cdots
\Theta_{n,n-2} \mP_{n,n-2}(1+\ga_{n,n-1})^*\\
&=& \Theta_{n,s} \mP_{n,s} \Theta_{n,s+1} \mP_{n,s+1}\cdots
\Theta_{n,n-1} \mP_{n,n-1},
\end{eqnarray*}
 by Corollary \ref{a7May9}.  $\Box $

$\noindent $

 Using Lemma \ref{a24May9}, we can strengthen
Theorem  \ref{A10May9}.

\begin{theorem}\label{B10May9}
\begin{enumerate}
\item $(1+\ga_n)^* = \Theta_{n,1}' \mE_{n,1} \Theta_{n,2}'
\mE_{n,2}\cdots \Theta_{n,n-1}' \mE_{n,n-1} $.  Moreover,  for
$s=1, \ldots , n-1$, $(1+\ga_{n,s})^* = \Theta_{n,s}' \mE_{n,s}
\Theta_{n,s+1}' \mE_{n,s+1}\cdots \Theta_{n,n-1}' \mE_{n,n-1} $.
\item $(1+\ga_n)^* =
 \Theta_{n,1}' \mP_{n,1} \Theta_{n,2}'
\mP_{n,2}\cdots \Theta_{n,n-1}' \mP_{n,n-1} $. Moreover,  for
$s=1, \ldots , n-1$, $(1+\ga_{n,s})^* =
 \Theta_{n,s}' \mP_{n,s} \Theta_{n,s+1}'
\mP_{n,s+1}\cdots \Theta_{n,n-1}' \mP_{n,n-1} $.
\end{enumerate}
\end{theorem}

{\it Proof}. The statements follow from Lemma \ref{a24May9},
(\ref{Tpns}),  and Theorem  \ref{A10May9}: repeat the proof of
Theorem \ref{A10May9}
 replacing $\Theta_{n,t}$ by $\Theta_{n,t}'$ everywhere for all
$t$. $\Box$

$\noindent $

By Theorem  \ref{A10May9}.(1) and Lemma \ref{a24May9}, the group
$\D_{n,s}$ is a free abelian group of rank ${n\choose s+1}s$  for
$s=1, \ldots , n-1$: 
\begin{equation}\label{Dns}
\D_{n,s}:=(1+\ga_{n,s})^*/\G_{n,s}=\Theta_{n,s}'\G_{n,s}/\G_{n,s}\simeq
\prod_{|J|=s+1}\prod_{j\in J\backslash \max (J)}\langle \th_{\max
(J), j}\rangle \simeq \Z^{{n\choose s+1}s},
\end{equation}
where the double product is the direct product of groups.

\begin{corollary}\label{a30May9}
$\CZ_{n,s}\simeq \mY_{n,s}'\simeq \Z^{{n\choose s+1}}$ for $s=1,
\ldots , n-1$ (see (\ref{apns1})).
\end{corollary}

{\it Proof}. Recall that
$\psi_{n,s}'((1+\ga_{n,s})^*)=\psi_{n,s}'(\Theta_{n,s})$ (Theorem
\ref{24May9}), $\mX_{n,s}'=\mK_{n,s}'\bigoplus \mY_{n,s}'$, and
$\psi_{n,s}'(\Theta_{n,s})=\mK_{n,s}'$, by (\ref{ptijJ1}). Then
\begin{eqnarray*}
 \CZ_{n,s}&=& \frac{\prod_{|I|=s}(1+\bgp_I)^*}{\psi_{n,s}((1+\ga_{n,s})^*)}\simeq
 \frac{\prod_{|I|=s}(1+\bgp_I)^*/\bG_{n,s}}{\psi_{n,s}((1+\ga_{n,s})^*)/\bG_{n,s}}\simeq  \frac{\mX_{n,s}'}{\psi_{n,s}'((1+\ga_{n,s})^*)}
 \simeq \frac{\mK_{n,s}'\bigoplus \mY_{n,s}'}{\psi_{n,s}'(\Theta_{n,s})}\\
 &=&\frac{\mK_{n,s}'\bigoplus \mY_{n,s}'}{\mK_{n,s}'}\simeq  \mY_{n,s}'\simeq \Z^{{n\choose
 s+1}}. \;\;\;\; \Box
\end{eqnarray*}

\begin{theorem}\label{30May9}
$\Upsilon_n= \G_n$ and $\Upsilon_{n,s}=\G_{n,s}$ for all $s=1,
\ldots , n$. In particular, the groups $\G_{n,s}$ are
$G_n$-invariant (hence, normal) subgroups of $\mS_n^*$ (since
$\Upsilon_{n,s}$ are so).
\end{theorem}

{\it Proof}. By Theorem  \ref{A10May9} and Lemma \ref{a24May9},
$(1+\ga_{n,s})^* = \Theta_{n,s}\G_{n,s} = \Theta_{n,s}'\G_{n,s}$
for  $s=1, \ldots , n$, and the last product is exact. Since
$\G_{n,s} \subseteq \Upsilon_{n,s}$, we have the equality
$(1+\ga_{n,s})^* =  \Theta_{n,s}'\Upsilon_{n,s}$. So, in order to
show that the equality $\G_{n,s} =\Upsilon_{n,s}$ holds it
suffices to prove that $\Theta_{n,s}'\cap \Upsilon_{n,s} = \{
1\}$. To prove this equality, first, we use an induction on $n\geq
2$, and then, for a fixed $n$, we  use a second downward induction
on $s=1, \ldots , n-1$, starting with $s=n-1$. For $n=2$, there is
a single option to consider, $(n,s) = (2,1)$. In this case, the
equality holds by Corollary \ref{a7May9}. Let $n>2$, and suppose
that the equality holds for all pairs $(n', s)$ with $n'<n$. For
$(n,n-1)$, the equality is true by Corollary \ref{a7May9}. Suppose
that $s<n-1$, and that the equality holds for all pairs $(n,s')$
with $s'=s+1, \ldots , n-1$. Suppose that $\Theta_{n,s}'\cap
\Upsilon_{n,s}\neq \{ 1\}$, we seek a contradiction. Choose an
element, say $u$, from the intersection such that $u\neq 1$. Then
the element $u$ is a unique product $u=\prod_{|J|=s+1}\prod_{j\in
J\backslash \max (J)}\th_{\max (J), j}(J)^{n(j,J)}$ where
$n(j,J)\in \Z$. Since $u\neq 1$, $n(j,J)\neq 0$ for some pair
$(j,J)$. Since $|J|=s+1<n$, the complement $CJ$ of the set $J$ is
a non-empty set. Let $f$ be the composition of the obvious algebra
homomorphisms: $$\mS_n\ra \mS_n/\sum_{i\in CJ}\gp_i\simeq \mS_J\t
L_{CJ}\ra \mS_J\t Q_{CJ}$$ where $Q_{CJ}$ is the field of
fractions of the Laurent polynomial algebra $L_{CJ}$. The algebra
$\mS_J\t Q_{CJ}$ is isomorphic to the algebra $\mS_{s+1}$ but over
the field $Q_{CJ}$. Let $\Theta_{s+1, s, J}'$ and $\G_{s+1, s, J}$
be the corresponding  $\Theta_{s+1, s}'$ and $\G_{s+1,s}$ for the
algebra $\mS_J\t Q_{CJ}\simeq \mS_{s+1}\t Q_{CJ}$ (over the field
$Q_{CJ}$). Since $f(\G_{n,s})\subseteq \G_{s+1, s,J}$,
$f(\Theta_{n,s}')\subseteq \Theta_{s+1, s, J}'\cdot U= \prod_{k\in
J\backslash \max (J)}\th_{\max (J), k}(J)^{n(k,J)}\cdot U$ where
$U:= \prod_{l\in J}U_{J\backslash l}(Q_{CJ})\subseteq \G_{s+1, s,
J}$, using the induction on $n$, the inclusion $f(u) \in
f(\Theta_{n,s}')\cap f(\G_{n,s})$ yields $n(j,J)=0$, a
contradiction. Therefore, $\Theta_{n,s}'\cap \Upsilon_{n,s} = \{
1\}$ and the statements of the theorem hold.  $\Box $

By (\ref{Dns}) and Theorem \ref{30May9}, 
\begin{equation}\label{Dns1}
(1+\ga_{n,s})^*/ \Upsilon_{n,s} = (1+\ga_{n,s})^*/ \G_{n,s} \simeq
\Z^{{n\choose s+1}s}.
\end{equation}
The next theorem gives explicit generators for the groups
$\mS_n^*$, $(1+\ga_n)^*$, and $(1+\ga_{n,s})^*$.
\begin{theorem}\label{25May9}
\begin{enumerate}
\item The group $(1+\ga_n)^*$ is generated by the following
elements:
\begin{enumerate}
\item $\th_{\max (J), j}(J)$ where $j\in J\backslash \max (J)$ and
$|J|=2, \ldots , n$; \item  $1+x_i^tE_{0\alpha}(I)$,
$1+x_i^tE_{\alpha 0}(I)$, $1+y_i^tE_{0\alpha}(I)$, and
$1+y_i^tE_{\alpha 0}(I)$ where $t\in \N\backslash \{ 0\}$,
$i\not\in I$, $|I|=1, \ldots , n-1$, $\alpha \in \N^I\backslash \{
0\}$; and \item $1+(\l -1)E_{00}(I)$, $1+E_{0\alpha} (I)$, and
$1+E_{\alpha 0}(I)$ where $\l \in K^*$, $I\neq \emptyset$, and
$\alpha \in \N^I\backslash \{ 0\}$.
\end{enumerate}
\item For $s=1, \ldots , n-1$, the group $(1+\ga_{n,s})^*$ is
generated by the following elements:
\begin{enumerate}
\item $\th_{\max (J), j}(J)$ where $j\in J\backslash \max (J)$ and
$|J|=s+1, \ldots , n$; \item  $1+x_i^tE_{0\alpha}(I)$,
$1+x_i^tE_{\alpha 0}(I)$, $1+y_i^tE_{0\alpha}(I)$, and
$1+y_i^tE_{\alpha 0}(I)$ where $t\in \N\backslash \{ 0\}$,
$i\not\in I$, $|I|=s, \ldots , n-1$, $\alpha \in \N^I\backslash \{
0\}$; and \item $1+(\l -1)E_{00}(I)$, $1+E_{0\alpha} (I)$, and
$1+E_{\alpha 0}(I)$ where $\l \in K^*$, $|I|=s, \ldots , n$, and
$\alpha \in \N^I\backslash \{ 0\}$.
\end{enumerate}
For $s=n$, the group $(1+\ga_{n,n})^*= (1+F_n)^*$ is generated by
the elements $1+(\l -1)E_{00}(I)$, $1+E_{0\alpha} (I)$, and
$1+E_{\alpha 0}(I)$ where $\l \in K^*$, $I=\{ 1, , \ldots , n\} $,
and $\alpha \in \N^n \backslash \{ 0\}$.
 \item The group $\mS_n^*= K^*\times (1+\ga_n)^*$
is generated by the elements from statement 1 and $K^*$.
\end{enumerate}
\end{theorem}

{\it Proof}. 1. Statement 1 is a particular case of statement 2
when $s=1$.

2. The statement is obvious for $s=n$. So, let $s=1, \ldots ,
n-1$. By Theorem \ref{B10May9}.(1), the group $(1+\ga_{n,s})^*$ is
generated by the sets $\Theta_{n,t}'$, $\mE_{n,t}$ where $t=1,
\ldots , n-1$. Each element of any of the sets $\Theta_{n,t}'$ is
a product of elements from (a). Recall that
$\mE_{n,t}:=\prod_{|I|=t}U_I(K)\ltimes E_\infty (\mS_{CI})$. Each
element of any of the groups $U_I(K)$ is a product of elements
from (c). For each $i=1, \ldots , n$, the algebra $\mS_1(i)$ is
the direct sum $\bigoplus_{j\geq 1}Ky_i^j\bigoplus K\bigoplus
\bigoplus_{j\geq 1} Kx_i^j\bigoplus F(i)$, see (\ref{mS1d}). By a
straightforward computation,
$$[1+aE_{\alpha\beta}(I), 1+bE_{\xi \rho}(I)]=\d_{\beta
\xi}(1+abE_{\alpha\beta}(I))$$ for all $a,b\in \mS_{CI}$ and
$\alpha, \beta , \xi, \rho \in \N^I$ where $[u,v]:=uvu^{-1}v^{-1}$
is the (group) commutator of elements $u$ and $v$. In this paper
the commutator stands for the group commutator (unless it is
stated otherwise).  For all $\l \in K^*$, $I$ with $|I|=s, \ldots
, n$, and $\alpha \in \N^I\backslash \{ 0\}$,
\begin{eqnarray*}
(1+(\l -1)E_{00}(I))\cdot (1+E_{0\alpha}(I))\cdot  (1+(\l -1)E_{00}(I))^{-1} &= & 1+\l E_{0\alpha}(I),  \\
(1+(\l -1)E_{00}(I))^{-1}\cdot (1+E_{\alpha 0}(I))\cdot (1+(\l
-1)E_{00}(I)) &= & 1+\l E_{\alpha 0}(I).
\end{eqnarray*}
It follows from these three facts that each element of any of the
sets $E_\infty (\mS_{CI})$ is a product of elements from sets (b)
and (c). The proof of statement 1 is complete.

3. Statement 3 is obvious.  $\Box $

The next theorem presents explicit generators for the group $G_n$.
\begin{theorem}\label{A25May9}
Let $J_s:= \{ 1, \ldots , s\}$ where $s=1, \ldots , n$. The group
$G_n= S_n\ltimes \mT^n\ltimes \Inn (\mS_n)$ is generated by the
transpositions $(ij)$ where $i<j$; the elements $t_{(\l , 1,
\ldots , 1)}:x_1\mapsto \l x_1$, $y_1\mapsto \l^{-1}y_1$,
$x_k\mapsto x_k$, $y_k\mapsto y_k$, $k=2, \ldots , n$; and the
inner automorphisms $\o_u$ where $u$ belongs to the following
sets:
\begin{enumerate}
\item $\th_{s,1}(J_s)$, $s=2, \ldots , n$; \item
$1+x_n^tE_{0\alpha}(J_s)$, $1+x_n^tE_{\alpha 0}(J_s)$,
$1+y_n^tE_{0\alpha}(J_s)$, and $1+y_n^tE_{\alpha 0}(J_s)$ where
$t\in \N\backslash \{ 0\}$, $s=1, \ldots , n-1$, and  $\alpha \in
\N^s\backslash \{ 0\}$; and \item  $1+(\l -1)E_{00}(J_s)$,
$1+E_{0\alpha} (J_s)$, and $1+E_{\alpha 0}(J_s)$ where $\l \in
K^*$,  $s=1, \ldots , n$, and $\alpha \in \N^s\backslash \{ 0\}$.
\end{enumerate}
\end{theorem}

{\it Proof}. The group $G_n= S_n\ltimes \mT^n\ltimes \Inn (\mS_n)$
(Theorem \ref{aInt24Apr9}.(3))  is generated by its three
subgroups: $S_n$, $\mT^n$, and $\Inn (\mS_n)=\{\o_v\, | \, v\in
(1+\ga_n)^*\}$. The transpositions generate the symmetric group
$S_n$. Then, by conjugating,
$$ (1i)t_{(\l , 1,
\ldots , 1)}(1i)^{-1} =t_{(1, \ldots , 1, \l , 1, \ldots , 1)}
\;\; (\l \; {\rm is \; on}\; i'{\rm th \; place})$$ we obtain
generators for the torus $\mT^n$.  Similarly, by conjugating the
elements of the sets 1, 2, and 3 (i.e. using $s \o_vs^{-1} = \o_{s
(v) }$ for all $s\in S_n$) we obtain all the elements from the
sets (a), (b) and (c) of Theorem \ref{25May9} when we identify the
groups $\Inn (\mS_n)$ and $(1+\ga_n)^*$ via $\o_v\lra v$. Now, the
theorem is obvious.
 $\Box $

$\noindent $


\section{The commutants of the groups $G_n$ and $\mS_n$,  and an analogue of the Jacobian homomorphism}\label{COMGJAC}

In this section, the groups $[G_n,G_n]$ and $G_n/[G_n,G_n]$ are
found (Theorem \ref{6Jun9}) and they are used to show the
uniqueness of an analogue $\mJ_n$ (see (\ref{SnJac})) of the
Jacobian homomorphism for $n>2$, and in finding the exotic
Jacobians $\mJ_n^{ex}$  for $n=1,2$.

{\bf The groups $[G_n,G_n]$ and $G_n/[G_n,G_n]$}.  The subgroup of
a group $G$ generated by all the commutators $[a,b]:=
aba^{-1}b^{-1}$ where $a,b\in G$ is called the {\em commutant} of
the group $G$ denoted either by $[G,G]$ or $G^{(1)}$. The
commutant is the least normal subgroup $G'$ of $G$ such that the
factor group $G/G'$ is abelian.
 If $\v :G\ra H$ is a
group homomorphism then $\v ([G,G])\subseteq [H,H]$. If, in
addition, the group $H$ is abelian then $[G,G]\subseteq \ker (\v
)$. To find the commutant of a group is a technical process
especially if the group is large. In the next two easy  lemmas we
collect patterns that appear in findings the commutant of the
group $G_n$. Their repeated applications make arguments short.

\begin{lemma}\label{a6Apr9}
\begin{enumerate}
\item The commutant $[A\ltimes B, A\ltimes B]$  of a skew product
$A\ltimes B$ of two groups is equal to $[A,A]\ltimes ([A,B]\cdot
[B,B])$ where $[A,B]$ is the subgroup of $B$ generated by all the
commutators $[a,b]:= aba^{-1}b$ for $a\in A$ and $b\in B$. Hence,
$B\cap [A\ltimes B, A\ltimes B]=[A,B]\cdot [B,B]$ and
$\frac{A\ltimes B}{[A\ltimes B, A\ltimes B]}\simeq
\frac{A}{[A,A]}\times \frac{B}{[A,B]\cdot [B,B]}$. \item
 If, in addition, the group $B$ is a direct product of groups
 $\prod_{i=1}^mB_i$ such that $aB_ia^{-1} \subseteq B_i$ for all
 elements $a\in A$ and $i=1, \ldots , m$. Then
 $[A\ltimes B, A\ltimes B]= [A,A]\ltimes \prod_{i=1}^m([A,B_i][B_i, B_i]).$
\end{enumerate}

\end{lemma}

{\it Proof}. 1.  Note that $[a,b]= \o_a(b) b^{-1}$ where $\o_a (b)
:= aba^{-1}$. For $a\in A$ and $b, c\in B$,
\begin{eqnarray*}
 c[a,b]&=&c\o_a(b)b^{-1} = \o_a( \o_{a^{-1}}(c) b) (\o_{a^{-1}}(c)b)^{-1} \o_{a^{-1}}(c) bb^{-1} \\
 &=&\o_a( \o_{a^{-1}}(c) b) (\o_{a^{-1}}(c)b)^{-1}\cdot
 \o_{a^{-1}}(c)\\
 &=& [ a,\o_{a^{-1}}(c) b]\cdot \o_{a^{-1}}(c).
\end{eqnarray*}
It follows from these equalities (when, in addition, we chose
$c\in [B,B]$) that the subgroup of $B$ which is generated by its
two subgroups, $[A,B]$ and $[B,B]$, is equal their set theoretic
product $[A,B][B,B]:=\{ ef\, | \, e\in [A,B], f\in [B,B]\}$. Then
the subgroup of $C:=[A\ltimes B, A\ltimes B]$ which is  generated
by its three subgroups $[A,A]$, $[A,B]$, and $[B,B]$ is equal to
the RHS, say $R$, of the equality of the lemma. It remains to
prove that $C\subseteq R$. This inclusion follows from the fact
that,  for all $a_1,a_2\in A$ and $b_1, b_2\in B$,
\begin{equation}\label{a1b1}
 [a_1b_1, a_2b_2]= \o_{a_1}([b_1, a_2]) \o_{a_1a_2}([b_1, b_2]) [
a_1, a_2] \o_{a_2}([a_1, b_2])
\end{equation}
 which follows from the equalities
$[ab, c] = \o_a([b,c])[a,c]$ and $[a,b]^{-1} = [ b,a]$:
\begin{eqnarray*}
 [a_1b_1, a_2b_2]&=& \o_{a_1}([b_1, a_2b_2]) [ a_1, a_2b_2]=
 ([a_2b_2, a_1]\o_{a_1}([a_2b_2, b_1]))^{-1}\\
 & =& (\o_{a_2}(b_2,
 a_1]) [ a_2, a_1] \o_{a_1} ( \o_{a_2}([b_2, b_1]) [a_2,
 b_1])^{-1}\\
 &=&  \o_{a_1}([b_1, a_2]) \o_{a_1a_2}([b_1, b_2])
[ a_1, a_2] \o_{a_2}([a_1, b_2]).
\end{eqnarray*}
2. By statement 1, it suffices to show that $[A, \prod_{i=1}^m
B_i]=\prod_{i=1}^m [ A,B_i]$. The general case follows easily from
the case when $m=2$ (by induction). The case $m=2$ follows from
(\ref{a1b1}) where we put $b_1=1$, $a_1\in A$, $a_2\in B_1$, and
$b_2\in B_2$. $\Box$

\begin{lemma}\label{a8Jun9}
\begin{enumerate}
\item Let $\v : G\ra H$ be a group epimorphism such that $\ker (\v
)\subseteq [G,G]$. Then  $[G,G]= \v^{-1}([H,H])$. \item Let $N$ be
a normal subgroup of a group such that $N\subseteq [G,G]$ and the
factor group $G/N$ is abelian. Then then $N=[G,G]$.
\end{enumerate}
\end{lemma}

{\it Proof}. 1. Since $\v$ is an epimorphism with $\ker (\v )
\subseteq [G,G]$, the inclusion $\v^{-1}([H,H])\subseteq [G,G]$ is
obvious. Then the  composition of the group homomorphisms
$G\stackrel{\v }{\ra}H\ra H/[H,H]$ and the fact that the group
$H/[H,H]$ is abelian yield the opposite inclusion
$\v^{-1}([H,H])\supseteq [G,G]$.

2. Applying statement 1 to the group epimorphism $\v : G\ra G/N$
we get statement 2: $[G,G]=\v^{-1} ([G/N,G/N])=\v^{-1} (e)=\ker
(\v )$. $\Box $

$\noindent $

For all transpositions $(ij)\in S_n$ and elements $t_{(\l_1,
\ldots , \l_n)}\in \mT^n$, 
\begin{equation}\label{comijt}
[(ij),t_{ (\l_1, \ldots , \l_n) } ]=t_{ (1,\ldots, 1,
\l_i^{-1}\l_j, 1,\ldots , 1, \l_j^{-1}\l_i, 1,\ldots , 1) }
\end{equation}
where the elements $\l_i^{-1}\l_j$ and $\l_j^{-1}\l_i$ are on
$i$th and $j$th place respectively.

\begin{lemma}\label{a5Jun9}
For each natural number $n\geq 2$, $[S_n\ltimes \mT^n, S_n\ltimes
\mT^n]= [S_n,S_n]\ltimes \mT^n_1$ where $\mT^n_1:= \{ t_{(\l_1,
\ldots , \l_n)}\in \mT^n \, | \, \prod_{i=1}^n \l_i =1\}$.
\end{lemma}

{\it Proof}. Let $R$ and $L$ be the RHS and the LHS of the
equality. By Theorem \ref{a6Apr9}.(1), (\ref{a1b1}), and
(\ref{comijt}), $R\supseteq L$. To prove the reverse inclusion
consider two group epimorphisms:
\begin{eqnarray*}
 \v : S_n\ltimes \mT^n &\ra & K^*, \;\; (\s , t_{(\l_1,
\ldots , \l_n)})\mapsto \prod_{i=1}^n\l_i,   \\
  \psi  : S_n\ltimes \mT^n &\ra & S_n\ltimes \mT^n/\mT^n\simeq
  S_n, \;\; (\s, t_\l ) \mapsto \s.
\end{eqnarray*}
Then $R\subseteq \ker (\v )= S_n\ltimes \mT^n_1$ and $R\subseteq
\psi^{-1}([S_n,S_n])= [S_n,S_n]\ltimes \mT^n$, hence $R\subseteq
(S_n\ltimes \mT^n_1)\cap ([S_n,S_n]\ltimes \mT^n)=L$, as required.
$\Box$

$\noindent $

Let $J$ be a subset of the set $\{ 1, \ldots , n\}$ that contains
at least two elements, let  $i$ and $j$ be two distinct elements
of the set $J$, and let $\l \in K^*$. By multiplying out, we see
that 
\begin{equation}\label{meJij1}
\mu_{J \backslash i} (y_i) e_{J \backslash j}= e_{J \backslash j}
\mu_{J \backslash i} (x_i) =e_{J \backslash j}-e_J,
\end{equation}
\begin{equation}\label{meJij2}
\mu_{J \backslash i} (y_i)e_J= e_J\mu_{J \backslash i} (x_i) =0,
\end{equation}
\begin{equation}\label{meJij3}
\mu_{J \backslash j} (x_jy_j) = 1-e_J,
\end{equation}
\begin{equation}\label{meJij4}
e_J\mu_{J \backslash j } (\l )= \mu_{J \backslash j} (\l ) e_J=\l
e_J.
\end{equation}
Note that (where $\l \in K^*$) 
\begin{equation}\label{tmJ}
[\th_{ij}(J), \mu_{J\backslash j}(\l )]=\mu_J(\l^{-1})
\end{equation}
since (by direct computations, consider the four cases as in
(\ref{curtij}))
$$[\th_{ij}(J), \mu_{J\backslash j}(\l )]*x^\alpha = \begin{cases}
\l^{-1}x^\alpha & \text{if } \forall k\in J:\alpha_k=0,\\
x^\alpha& \text{otherwise}.\\
\end{cases} $$
Alternatively, using the equalities (\ref{meJij1}),
(\ref{meJij2}), (\ref{meJij3}), and (\ref{meJij4}), we can show
directly  that (\ref{tmJ}) holds:
\begin{eqnarray*}
 [\th_{ij}(J), \mu_{J\backslash j}(\l )]&=& \th_{ij}(J)\mu_{J\backslash j}(\l ) \th_{ji}(J)\mu_{J\backslash j}(\l^{-1})\\
  & =&
 \mu_{J\backslash i}(y_i) \cdot \mu_{J\backslash j}(x_j) \mu_{J\backslash j}( \l)\mu_{J\backslash j}(y_j)\cdot \mu_{J\backslash i}(x_i)
 \mu_{J\backslash j}(\l^{-1})\\
 & =&\mu_{J\backslash i}(y_i) \cdot \mu_{J\backslash j}(x_jy_j)\cdot  \mu_{J\backslash j}( \l)
\cdot \mu_{J\backslash i}(x_i)
 \mu_{J\backslash j}(\l^{-1}) \\
 &=& \mu_{J\backslash i}(y_i) \cdot (1-e_J)\cdot  \mu_{J\backslash j}( \l)
\cdot \mu_{J\backslash i}(x_i)
 \mu_{J\backslash j}(\l^{-1}) \;\;\; \;\;\; ({\rm by}\;\; (\ref{meJij3})) \\
 &=& (1+(\l -1)\mu_{J\backslash i}(y_i) e_{J\backslash j}  \mu_{J\backslash i}(x_i))
 \cdot  \mu_{J\backslash j}(\l^{-1}) \;\;\; \;\;\;\;\;\;  ({\rm by}\;\; (\ref{meJij2})) \\
 &=& (1+(\l -1) (e_{J\backslash j}-e_J) \mu_{J\backslash i}(x_i))
 \cdot  \mu_{J\backslash j}(\l^{-1}) \;\;\; \;\;\; \;\;\; \;\; ({\rm by}\;\; (\ref{meJij1})) \\
 &=& (1+(\l -1) (e_{J\backslash j}-e_J))\cdot  \mu_{J\backslash j}(\l^{-1}) \;\;\; \;\;\;\; \;\;\; \;\;\; \;\;\;({\rm by}\;\; (\ref{meJij1}),  (\ref{meJij2})) \\
 & = & (\mu_{J\backslash j}(\l
 )+(1-\l ) e_J)\cdot \mu_{J\backslash j}(\l^{-1}) = 1+(1-\l )
 \l^{-1}e_J \;\;\; \;\;\; ({\rm by}\;\;
 (\ref{meJij4}))\\
 &=& 1+(\l^{-1}-1)e_J= \mu_J(\l^{-1}).
\end{eqnarray*}
By taking the inverse of both sides of (\ref{tmJ})  and using the
fact that $[a,b]^{-1}=[b,a]$, we have the equality
\begin{equation}\label{tmJ1}
[\mu_{J\backslash j}(\l ), \th_{ij}(J) ]=\mu_J(\l ).
\end{equation}
Let $J$ be a subset of the set $\{ 1,\ldots , n\}$. If $i$ and $j$
are distinct elements of the set $J$ (hence $|J|\geq 2$)  then,
for all elements $s\in S_n$, 
\begin{equation}\label{sost}
s\o_{\th_{ij}(J)}s^{-1}=\o_{\th_{s(i)s(j)}(s(J))},
\end{equation}
\begin{equation}\label{sost1}
[(ij), \o_{\th_{ij}(J)}]=\o_{\th_{ij}(J)^{-2}}.
\end{equation}
The equality (\ref{sost}) is obvious, the equality  (\ref{sost1})
follows from (\ref{sost}) and (\ref{tiji}): $$ [(ij),
\o_{\th_{ij}(J)}]=
 (ij)\o_{\th_{ij}(J)}(ij)^{-1}\o_{\th_{ij}(J)}^{-1}=
 \o_{\th_{ji}(J)}\o_{\th_{ij}(J)^{-1}}=\o_{\th_{ij}(J)^{-1}}\o_{\th_{ij}(J)^{-1}}=
 \o_{\th_{ij}(J)^{-2}}.$$
If $i$, $j$, and $k$ are distinct elements of the set $J$ (hence
$|J|\geq 3$)  then 
\begin{equation}\label{sost2}
[(ik), \o_{\th_{ij}(J)}]=\o_{\th_{ki}(J)}.
\end{equation}
In more detail,
\begin{eqnarray*}
 [(ik), \o_{\th_{ij}(J)}]&=& (ik) \o_{\th_{ij}(J)}(ik)^{-1}
 \o_{\th_{ij}(J)}^{-1}=\o_{\th_{kj}(J)}\o_{\th_{ji}(J)}\;\;\;
 ({\rm by}\;\; (\ref{sost}), (\ref{tiji}))\\
&=&\o_{\th_{ki}(J)} \;\;\; \;\;\; \;\;\; ({\rm by}\;\; (\ref{tijjk})).\\
\end{eqnarray*}
By (\ref{sost2}), if $n>2$ then the current group $\Theta_n$
belongs to the commutant $[G_n,G_n]$, but for $n=2$ this is not
true (Theorem \ref{6Jun9}.(1)), and this is the reason for
existence of the exotic `Jacobian' homomorphism $\mJ_2^{ex}$.

Let $\th_{ij}:=\th_{ij}(\{ i,j\})$ and $\mu_j(\l ) :=\mu_{\{
j\}}(\l)$ where $\l\in K^*$. Then 
\begin{equation}\label{sost3}
[t_{(1, \ldots , 1,\l_i,1,\ldots , 1)},
\o_{\th_{ij}}]=\o_{\mu_j(\l_i^{-1})},
\end{equation}
where the scalar $\l_i\in K^*$ is on the $i$th place. In more
detail, $ [t_{(1, \ldots , 1,\l_i,1,\ldots , 1)}, \o_{\th_{ij}}] =
\o_{\mu_j(\l_i^{-1}y_i) \mu_i(x_j)}\cdot
\o_{\th_{ij}^{-1}}=\o_{\mu_j(\l_i^{-1}) \th_{ij}\th_{ij}^{-1}}=
\o_{\mu_j(\l_i^{-1})}$.


\begin{theorem}\label{6Jun9}
Let $\th := \th_{12}(\{ 1,2\})$ and $\CN_2:= \{ \o_u\, | \, u\in
\langle\th^2 \rangle \cdot  \prod_{|I|=1}U_I(K)\ltimes
E_I(\mS_{CI})\cdot (1+\ga_{2,2})^*\}\subseteq G_2$. Then
\begin{enumerate}
\item $[G_n,G_n]= \begin{cases}
\{ \o_u\, | \, u\in E_\infty (K)\} & \text{if }n=1,\\
\mT^1_1\ltimes  \CN_2& \text{if }n=2, \\
[S_n,S_n]\ltimes \mT^n_1\ltimes \Inn (\mS_n) & \text{if }n>2.\\
\end{cases}$
\item $ G_n/[G_n,G_n]\simeq \begin{cases}
K^*\times K^*& \text{if } n=1,\\
\Z/ 2\Z \times  K^* \times \Z/ 2\Z  &\text{if } n=2, \\
\Z/ 2\Z \times  K^*  &   \text{if } n>2. \\
\end{cases} $
\end{enumerate}
\end{theorem}

{\it Proof}.  Recall that $(1+\ga_n)^*\simeq \Inn (\mS_n)$, $u\lra
\o_u$ (Theorem \ref{aInt24Apr9}.(3)). To save on notation we
identify these two groups. Then $G_n=S_n\ltimes \mT^n\ltimes \Inn
(\mS_n) = S_n\ltimes \mT^n\ltimes (1+\ga_n)^*$.

The case $n=1$. By Theorem 4.1, \cite{shrekaut},
$$G_1\simeq \mT^1\ltimes (1+F)^*\simeq \mT^1\ltimes (U(K)\ltimes
E_\infty (K))= (\mT^1 \times U(K))\ltimes E_\infty (K).$$ Since
 $[E_\infty (K), E_\infty (K)]=E_\infty (K)$ (hence $E_\infty
(K)\subseteq [G_1, G_1]$) and the factor group $G_1/E_\infty (K)
\simeq \mT^1 \times U(K)$ is abelian, by Lemma \ref{a8Jun9}.(2),
$[G_1, G_1]=\{ \o_u\, | \, u\in E_\infty (K)\}$. Hence
$G_1/[G_1,G_1]\simeq \mT^1\times U(K)\simeq K^*\times K^*$.

Let $n\geq 2$. Note that $E_\infty (\mS_{CI})=[E_\infty
(\mS_{CI}), E_\infty (\mS_{CI})]\subseteq [G_n,G_n]$ for all
nonempty subsets $I$ of the set $\{ 1,\ldots , n\}$. It follows
from (\ref{tmJ1}), (\ref{sost2}), and Theorem \ref{B10May9}.(1)
that
$$ (1+\ga_{n,2})^* \subseteq [G_n,G_n].$$

The case $n=2$. By (\ref{sost1}), $\th^2\in[ G_2,G_2]$.  By
(\ref{sost3}),  $\prod_{|I|=1}U_I(K)\ltimes E_I(\mS_{CI})\subseteq
[ G_2,G_2]$. By Theorem \ref{B10May9}.(1) and (\ref{apns1}),
$$\overline{G}_2:=G_2/(1+\ga_{2,2})^*\simeq S_2\ltimes \mT^2\ltimes \langle
\th \rangle \ltimes \prod_{|I|=1}U_I(K)\ltimes E_I(L_{CI}).$$ Note
that $[(12), \o_\th ]=\o_{\th^{-2}}$ (by (\ref{sost1})), and,  for
all elements $t_\l \in \mT^2$,
$$[t_\l , \o_\th ]\equiv \o_{\mu_2(\l_1^{-1})\mu_1(\l_2)}\mod
(1+\ga_{2,2})^*.$$ Indeed,
\begin{eqnarray*}
 [t_\l , \o_\th ]&\equiv & \o_{\mu_2(\l_1^{-1})\th
 \mu_1(\l_2)}\o_{\th^{-1}}\equiv \o_{\mu_2(\l_1^{-1})\th
 \mu_1(\l_2)\th^{-1}}\equiv \o_{\mu_2(\l_1^{-1})
 \mu_1(\l_2)\th \th^{-1}}\\
&\equiv & \o_{\mu_2(\l_1^{-1})\mu_1(\l_2)}\mod (1+\ga_{2,2})^*.
\end{eqnarray*}
It follow that the group $N:=\langle \th^2\rangle \ltimes
\prod_{|I|=1}U_I(K)\ltimes E_I(L_{CI})$ is a normal subgroup of
$\overline{G}_2$, $N\subseteq [ \overline{G}_2,\overline{G}_2]$,
and $\overline{G}_2/N\simeq (S_2\ltimes \mT^2)\times \frac{\langle
\th \rangle }{\langle \th^2 \rangle}$. By Lemma \ref{a5Jun9},
$[\overline{G}_2/N, \overline{G}_2/N]= [ S_2\ltimes \mT^2,
S_2\ltimes \mT^2] = \mT^2_1$. Then, by Lemma \ref{a8Jun9}.(1),
statement 1 follows. Then, by Lemma \ref{a8Jun9}.(1),
$$\frac{G_2}{[G_2, G_2]}\simeq \frac{\overline{G}_2/N}{[\overline{G}_2/N,
\overline{G}_2/N]}\simeq S_2\times \frac{\mT^2}{\mT^2_1}\times
\frac{\langle \th \rangle }{\langle \th^2 \rangle}\simeq
\Z_2\times K^*\ltimes \Z_2. $$
 The case $n>2$. By Theorem \ref{B10May9}.(1), $(1+\ga_{n,1})^* =\Theta_{n,1}' \mE_{n,1}
\cdots \Theta_{n,n-1}' \mE_{n,n-1} $. By (\ref{sost2}),
$\Theta_{n,s}'\subseteq [ G_n, G_n]$ for all $s=1, \ldots , n-1$.
By (\ref{tmJ1}), $\mE_{n,s}\subseteq [G_n,G_n]$ for all $s=2,
\ldots , n-1$, and, by (\ref{sost3}), $\mE_{n,1}\subseteq [
G_n,G_n]$. Therefore, $(1+\ga_{n,1})^*\subseteq [G_n,G_n]$. Then
the factor group $\overline{G}_n:= G_n/(1+\ga_{n,1})^*$ is
isomorphic to the group $S_n\ltimes \mT^n$. By Lemma \ref{a5Jun9},
$[\overline{G}_n,\overline{G}_n]=[S_n,S_n]\ltimes \mT^n_1$, and
statement 1 follows, by Lemma \ref{a8Jun9}.(1). By Lemma
\ref{a8Jun9}.(1),
$$\frac{G_n}{[G_n, G_n]}\simeq \frac{\overline{G}_n}{[\overline{G}_n,
\overline{G}_n]}\simeq \frac{S_n}{[S_n,S_n]}\times
\frac{\mT^n}{\mT^n_1}\simeq \Z_2\times K^*. $$ The proof of the
 theorem is complete. $\Box$

Recall that $\aff_n:= S_n\ltimes \mT^n$.
\begin{corollary}\label{a12Jun9}
\begin{enumerate}
\item $\frac{\Inn (\mS_n)}{\Inn (\mS_n)\cap [G_n,G_n]}\simeq
\begin{cases}
K^*& \text{if }n=1,\\
\Z_2& \text{if }n=2,\\
0& \text{if }n>2.\\
\end{cases}$
\item $\frac{\aff_n}{[\aff_n, \aff_n]}\simeq \begin{cases}
K^*& \text{if }n=1,\\
\Z_2\times K^*& \text{if }n>1.\\
\end{cases}$
\item $\frac{G_n}{[G_n, G_n]}\simeq \frac{\aff_n}{[\aff_n,
\aff_n]}\times \frac{\Inn (\mS_n)}{\Inn (\mS_n)\cap
[G_n,G_n]}\simeq \frac{\aff_n}{[\aff_n, \aff_n]}\times
\begin{cases}
K^*& \text{if }n=1,\\
\Z_2& \text{if }n=2,\\
0& \text{if }n>2.\\
\end{cases}$
\end{enumerate}
\end{corollary}

{\it Proof}. 1. We keep the notation of the proof of Theorem
\ref{6Jun9} (in particular, we identify the groups $\Inn (\mS_n)$
and $(1+\ga_n)^*$). For $n=1$, $\Inn (\mS_1)=U(K)\ltimes E_\infty
(K)$ and $[G_1,G_1]=E_\infty (K)$ (Theorem \ref{6Jun9}.(1)) and
statement follows.

For $n=2$, by Theorem \ref{6Jun9}.(1), $\Inn (\mS_n)/\Inn
(\mS_n)\cap [G_n,G_n]\simeq \langle \th \rangle /\langle \th^2
\rangle \simeq \Z_2$.

For $n=3$, by Theorem \ref{6Jun9}.(1), $\Inn (\mS_n)\subseteq
[G_n,G_n]$.

2. For $n=1$, $\aff_1=\mT^1$ and statement 2 is obvious. For
$n>1$, statement 2 follows from Lemma \ref{a5Jun9}:
$\frac{\aff_n}{[\aff_n, \aff_n]}\simeq \frac{S_n}{[S_n,
S_n]}\times \frac{\mT^n}{\mT^n_1}\simeq \Z_2\times K^*$.

3. Since $G_n=\aff_n\ltimes \Inn (\mS_n)$, the first equality
follows from Lemma \ref{a6Apr9}.(1), and  then the second equality
follows from statement 1. $\Box $

$\noindent $

{\bf An analogue of the polynomial Jacobian homomorphism}. We keep
the notation of the Introduction.  We want to find an analogue of
the polynomial Jacobian homomorphism (\ref{polJac}) for the
algebra $\mS_n$. The algebra $\mS_n$ is noncommutative,
non-Noetherian, with trivial centre, i.e. $Z(\mS_n)=K$ Proposition
4.1, \cite{shrekalg}, and there are no obvious `partial'
derivatives for the algebra $\mS_n$. So, in order to find the
analogue we, first, define the Jacobian homomorphism in invariant
group-theoretic  terms, i.e. we select natural properties/
conditions that uniquely determine $\CJ$. Then, for the algebra
$\mS_n$,  the conditions obtained determine uniquely an analogue
of the Jacobian homomorphism  but for $n=1,2$ where there are
exactly two of them.

The group $\CP_n = \Aff_n\times_{ex} \S_n$ is an exact product of
its two subgroups where $\Aff_n:=\{ \s_{A,a}:x\mapsto Ax+a\, | \,
A\in \GL_n(K), a\in K^n\}$ is the {\em affine group} and
$$\S_n:=\{ \s \in \CP_n\, | \, \s (x_i)\equiv x_i \mod (x_1,
\ldots , x_n)^2, i=1, \ldots , n\}$$ is the {\em Jacobian} group
where $(x_1, \ldots , x_n) $ is the maximal ideal of the
polynomial algebra $P_n$. Recall that an exact product means that
each element $\s \in \CP_n$ is a unique product $\s =
\s_{A,a}\cdot \xi$ where $\s_{A,a}\in \Aff_n$ and $\xi \in \S_n$.
The Jacobian homomorphism $\CJ_n$ is determined by its restriction
to the affine subgroup since $\CJ_n (\xi ) = 1$ for all $\xi \in
\S_n$ (trivial); and  
\begin{equation}\label{JsdetA}
\CJ_n(\s )=\CJ_n (\s_{A,a})= \det (A).
\end{equation}
The affine group $\Aff_n$ is the semidirect product $U_n\ltimes
\SAff_n$ where $U_n:=\{ \tau_\l : x_1\mapsto \l x_1, x_i\mapsto
x_i, i>1\, | \, \l \in K^*\}$ and  $\SAff_n := \{ \s_{A,a}\, | \,
A\in \SL_n(K), a\in K^n\}$, and $\SAff_n= [ \Aff_n, \Aff_n]$ is
the commutant of the group $\Aff_n$. Then the commutant $[\CP_n,
\CP_n]$ of the group $\CP_n$ is $\SAff_n\times_{ex} \S_n$ (this is
trivial, since $[U_n, \S_n] = \S_n$), and $\CP_n = U_n\ltimes
[\CP_n, \CP_n]= (U_n\ltimes [ \Aff_n, \Aff_n])\times_{ex} \S_n$.
The Jacobian homomorphism is uniquely determined by the following
commutative diagram:
$$
\xymatrix{ \CP_n\ar[r]^{\CJ}\ar[d]  & K^* \\
 \frac{\CP_n}{[\CP_n, \CP_n]}\ar[r]^{\simeq}   & U_n\ar[u]^{\det} }
$$
The group $G_n=S_n\ltimes \mT^n\ltimes \Inn (\mS_n)$ has a similar
structure as the group $\CP_n$. The subgroup $\aff_n:=S_n\ltimes
\mT^n$ is an affine part of the group $G_n$ and the subgroup $\Inn
(\mS_n)$ plays a role of the Jacobian subgroup $\S_n$  since
\begin{itemize}
\item (Corollary 5.5, \cite{shrekaut}) $\Inn (\mS_n) = \{ \s \in
G_n\, | \, \s (x_i) \equiv x_i\mod \gp_i,\; \s (y_i) \equiv
y_i\mod \gp_i, \; \forall i\}.$
\end{itemize}

{\it Definition}. An analogue $\mJ_n$ of the polynomial Jacobian
homomorphism $\CJ_n$ is a group homomorphism $\mJ_n:G_n\ra K^*$
which acts on the affine subgroup $\aff_n$ as in the polynomial
case.

$\noindent $

There is at least one such a map which is given in (\ref{SnJac})
and (\ref{SnJac1}).

\begin{theorem}\label{15Jun9}
\begin{enumerate}
\item For $n>2$, the analogue $\mJ_n$ of the polynomial Jacobian
homomorphism $\CJ_n$ is unique and given in  (\ref{SnJac}) and
(\ref{SnJac1}). \item For $n=1,2$, there is another one
$\mJ_n^{ex}$, the so-called exotic Jacobian homomorphism, given by
the rule
\begin{enumerate}
\item For $n=1$, $\s = t_\l \cdot \o_u\in G_1=\mT^1\ltimes \{
\o_v\, | \, v\in (1+F)^*\}$ where $t_\l \in \mT^1$ and $u\in
(1+F)^*\simeq \GL_\infty (K)$, $\mJ_1^{ex}(\s ) = \l \cdot \det
(u)$. The homomorphisms $\mJ_1$ and $\mJ_1^{ex}$ are algebraically
independent characters of the group $G_1$. \item For $n=2$, $\s =
st_\l \o_{\th^i}\xi \in G_2$ where $s\in S_2$, $t_\l \in \mT^2$,
 $i\in \{ 0,1\}$, and $\xi \in \CN_2$, $\mJ_2^{ex} (\s ) =
 (-1)^i\sign (s) \l_1\l_2$. Note that $(\mJ_2^{ex})^2= \mJ_2^2$.
\end{enumerate}
\end{enumerate}
\end{theorem}

{\it Proof}. 1. Statement 1 follows from the fact that
$\frac{G_n}{[G_n, G_n]}\simeq \frac{\aff_n}{[\aff_n, \aff_n]}$
(Corollary \ref{a12Jun9}.(3)).

2. Statement 2 follows from Corollary \ref{a12Jun9}.(3)  $\Box $

$\noindent $

Department of Pure Mathematics

University of Sheffield

Hicks Building

Sheffield S3 7RH

UK

email: v.bavula@sheffield.ac.uk


\begin{thebibliography}{99}

\bibitem{Bass-book-K-theory} H. Bass,  Algebraic $K$-theory. W. A. Benjamin, Inc., New York-Amsterdam,
 1968, 762 pp.




\bibitem{Bav-Jacalg} V. V.  Bavula, The Jacobian algebras, {\em J. Pure Appl.
Algebra}, {\bf 213} (2009),  664-685;  Arxiv:math.RA/0704.3850.




\bibitem{A1rescen} V. V. Bavula, The group of automorphisms of the first Weyl algebra in prime
characteristic and the restriction map,  {\em Glasgow Math. J.},
{\bf 5} (2009), 263-274;   Arxiv:math.RA/0708.1620.



\bibitem{shrekalg} V. V. Bavula, The algebra of one-sided inverses of a polynomial algebra,  ArXiv:math.AG/0903.0641.

\bibitem{shrekaut} V. V. Bavula, The group of automorphisms of the
algebra of one-sided inverses of a polynomial algebra,
ArXiv:math.AG/0903.3049.

\bibitem{jacaut} V. V. Bavula, The group of automorphisms of the Jacobian algebra
$\mA_n$.


\bibitem{K1aut} V. V. Bavula, ${\rm K}_1(\mS_1)$ and the group of automorphisms of the  algebra
$\mS_2$  of one-sided inverses of a polynomial algebra in two
variables, arXiv:0906.0600.



\bibitem{Dix}
J. Dixmier,  Sur les alg\`{e}bres de Weyl, {\it Bull. Soc. Math.
France} {\bf 96} (1968), 209--242.

\bibitem{Jacobson-StrRing} N. Jacobson, ``Structure of rings,''Am.
Math. Soc. Colloq., Vol. XXXVII (rev. ed.), Am. Math. Soc.,
Providence, 1968.




\bibitem{Mak-LimBSMF84} L. Makar-Limanov, On automorphisms of Weyl algebra, {\em Bull. Soc. Math.
France}, {\bf 112} (1984), 359--363.

\bibitem{Milnor-book-K-theory}J. Milnor,  Introduction to algebraic $K$-theory. Annals of Mathematics Studies,
 No. 72. Princeton University Press, Princeton, N.J.; University of Tokyo Press, Tokyo, 1971, 184 pp.

\bibitem{jung} H. W. E. Jung, Uber ganze birationale
Transformationen der Ebene,  {\it J. Reine Angew. Math.}, {\bf
184} (1942), 161-174.

\bibitem{kulk} W. Van der Kulk, On polynomial rings in two
variables, {\it Nieuw Arch. Wisk. } (3) {\bf 1} (1953), 33-41.

\bibitem{Swan} R. G. Swan,  Algebraic $K$-theory. Lecture Notes in Mathematics,
 No. 76 Springer-Verlag, Berlin-New York,  1968, 262 pp.

\end{thebibliography}
\end{document}